\documentclass[11pt,reqno]{amsart}

\usepackage[T2A]{fontenc}
\usepackage[utf8]{inputenc}
\usepackage[russian,english]{babel}

\usepackage{amsthm,amsfonts,amssymb,amsmath,oldgerm,mathrsfs,mathabx,enumitem}
\usepackage[scr=boondoxo]{mathalfa}
\numberwithin{equation}{section}
\usepackage{fullpage,setspace,fancyhdr,bm}
\usepackage{graphicx,psfrag,listings} 
\usepackage[pdftex,dvipsnames]{xcolor}
\usepackage{cancel}

\usepackage[top=30mm,bottom=30mm,left=25mm,right=25mm,a4paper]{geometry}

\usepackage{hyperref}
\usepackage{graphics}

\setlist[itemize,1]{label=\ensuremath{\diamond}}

\theoremstyle{plain}
\newtheorem{theorem}{Theorem}

\newtheorem{proposition}[theorem]{Proposition}
\newtheorem{lemma}[theorem]{Lemma}

\theoremstyle{definition}

\newtheorem{remark}[theorem]{Remark}

\newcommand{\de}{\mathrm{d}}

\newcommand\dd{\mathrm{d}}

\newcommand{\re}{\mathrm{e}}
\newcommand{\ri}{\mathrm{i}}

\DeclareMathOperator{\Real}{Re}
\DeclareMathOperator{\Imag}{Im}

\newcommand{\norm}[1]{\lVert#1\rVert}

\newcommand{\absolute}[1]{\lvert#1\rvert}
\newcommand{\Absolute}[1]{\left\lvert#1\right\rvert}
 
\newcommand{\Scalp}[1]{\left\langle #1\right\rangle}

\newcommand{\transp}[1]{#1^\mathrm{T}}
\newcommand{\Set}[1]{\left\lbrace#1\right\rbrace}

\newcommand{\CC}{{\mathbb C}}

\newcommand{\RR}{{\mathbb R}}
\renewcommand{\SS}{{\mathbb S}}

\newcommand\cC{{\mathcal C}}

\newcommand{\eu}{\re}
\newcommand{\iu}{\ri}

\title{Orbital Stability of Plane Waves in the Klein-Gordon Equation against Localized Perturbations}

\author{Emile Bukieda}
\email{{\tt emile.bukieda@kit.edu}}

\author{Louis Gar\'enaux}
\email{{\tt louis.garenaux@inria.fr}}

\author{Bj\"orn de Rijk}
\email{{\tt bjoern.rijk@kit.edu}}

\begin{document}

\begin{abstract}
We investigate the stability and long-term behavior of spatially periodic plane waves in the complex Klein-Gordon equation under localized perturbations. Such perturbations render the wave neither localized nor periodic, placing its stability analysis outside the scope of the classical orbital stability theory for Hamiltonian systems developed by Grillakis, Shatah, and Strauss. Inspired by Zhidkov's work on the stability of time-periodic, spatially homogeneous states in the nonlinear Schr\"odinger equation, we develop an alternative method that relies on an amplitude-phase decomposition and leverages conserved quantities tailored to the perturbation equation. We establish an orbital stability result of plane waves that is locally uniform in space, accommodating $L^2$-localized perturbations as well as unbounded phase modulations. Our result is sharp in the sense that it holds up to the spectral stability boundary. 

\vspace{0.5em}


\vspace{0.5em}

\end{abstract}

\maketitle

\section{Introduction}

This paper focuses on the stability and modulational behavior of plane-wave solutions to the one-dimensional complex Klein-Gordon equation
\begin{align} \label{e:KG}
u_{tt} - u_{xx} + f\big(|u|^2\big) u = 0, \qquad x, t \in \RR, \quad u(x,t) \in \mathbb{C},
\end{align}
with nonlinearity $f \in \cC^2([0,\infty),\RR)$. Plane waves represent the most basic periodic traveling-wave solutions. They are monochromatic waves, characterized by a single nonzero Fourier mode, of the form
\begin{align} \label{e:plane-wave}
u(x,t) = a \re^{\ri k x + \ri \omega t},
\end{align}
where $a > 0$ denotes the amplitude, $k \in \RR \setminus \{0\}$ is the wave number, $\omega \in \RR$ represents the temporal frequency, and $s = -\smash{\tfrac{\omega}{k}}$ is the associated wave speed. Inserting~\eqref{e:plane-wave} into~\eqref{e:KG}, we find that~\eqref{e:plane-wave} is a solution to~\eqref{e:KG} if and only if 
\begin{align} \label{e:disp-rel}
\omega^2 = k^2 + f(a^2).
\end{align}
The nonlinear dispersion relation~\eqref{e:disp-rel} expresses the temporal frequency $\omega$ in terms of the wave number $k$ and the amplitude $a$. 

The Klein-Gordon equation arises in a wide array of physical contexts, including nonlinear wave propagation, superconductivity, and quantum field dynamics; see, for instance, the classical textbooks~\cite{Lee1981Particle,Reed1975Methods,Whitham1999Linear} and references therein. In many applications, the nonlinearity takes a specific form, such as the power law $f(\nu) =\smash{1 \pm \nu^p}$ with $p \in \mathbb{N}$. In this work, however, we consider a general class of nonlinearities and impose no specific structural assumptions on $f$, beyond its $\cC^2$-regularity. 

Since the plane-wave solution~\eqref{e:plane-wave} traces a circle of radius $a$ centered at the origin, only the leading-order behavior of the nonlinearity near this circle is relevant for our stability analysis. This behavior is naturally described in polar coordinates by writing $u = a\re^{\rho + \ri \phi}$, where $\rho$ measures deviations from the circle. Expanding the nonlinearity yields
\begin{align} \label{e:fexpans}
f\big(|u|^2\big) u = f\big(a^2\big) u + 2a^2 f'\big(a^2\big) \rho u + \mathcal{O}\big(\rho^2 u\big),
\end{align}
so that for small $\rho$, the dominant contributions are the linear term $f(a^2) u$ and the quadratic term $2a^2f'(a^2)\rho u$. In the context of quantum field theory, the linear term can be interpreted as a mass term: $f(a^2)>0$ corresponds to a real mass, $f(a^2)=0$ to a massless field, and $f(a^2)<0$ to an imaginary, or \emph{tachyonic}, mass. Moreover, the sign of $f'(a^2)$ determines the nature of the leading-order nonlinear interactions: it is \emph{defocusing} if $f'(a^2)>0$ and \emph{focusing} if $f'(a^2)<0$. As we will see, the signs of $f(a^2)$ and $f'(a^2)$ are decisive for the stability properties of the plane wave and therefore play a central role in our analysis. 

The Klein-Gordon equation is a prototypical example of a Hamiltonian system
\begin{align} \label{e:HamSys}
\mathbf{u}_t = J H'(\mathbf{u}),
\end{align}
where $J$ is a skew-symmetric operator and $H$ is a nonlinear functional on a suitable Hilbert space. Upon introducing the vector $\mathbf{u} = (\Real(u), \Imag(u),\Real(u_t), \Imag(u_t))^\top$, the Klein-Gordon equation~\eqref{e:KG} can indeed be cast into the form~\eqref{e:HamSys} with
\begin{align*}
J = \begin{pmatrix} 0 & I_2 \\ -I_2 & 0\end{pmatrix}, \qquad I_2 = \begin{pmatrix} 1 & 0 \\ 0 & 1 \end{pmatrix},
\end{align*}
and Hamiltonian
\begin{align*}
H(\xi,\zeta,\mu,\nu) = \frac{1}{2} \int_{\mathcal{I}} \xi'(x)^2 + \zeta'(x)^2 + \mu(x)^2 + \nu(x)^2 + F\big(\xi(x)^2 + \zeta(x)^2 \, \de x,
\end{align*}
where $F$ is the primitive function of $f$ with $F(0) = 0$. Here, we have $\mathcal{I} = \RR$ in case of localized initial data, and $\mathcal{I} = [0,L]$ in case of $L$-periodic initial data. 

A characteristic feature of Hamiltonian systems is the conservation of the Hamiltonian along solutions. In addition to this, the complex Klein-Gordon equation exhibits a gauge (or rotational) invariance, i.e., if $u(t)$ is a solution to~\eqref{e:KG}, then so is $\re^{\ri \phi} u(t)$ for any $\phi \in \RR$. This gives rise to a second conserved quantity: the \emph{charge}, defined by
\begin{align*}
Q(\xi,\zeta,\mu,\nu) = \int_{\mathcal{I}} \xi(x)\nu(x) - \zeta(x)\mu(x) \, \de x.
\end{align*}

Conservation laws play a crucial role in understanding the long-term dynamics of Hamiltonian systems. Their use in proving global existence or finite-time blow-up of solutions is classical; see, for example,~\cite{cazenave1998introduction} for applications to the Klein–Gordon equation. Conserved quantities also form the backbone of the orbital stability theory of traveling (or standing) wave solutions. The primary approach to proving orbital stability of traveling waves is to characterize them as constrained minimizers of the Hamiltonian. This variational approach, dating back to ideas of Boussinesq, was first introduced in the nonlinear stability analysis of solitary waves in the Korteweg-de Vries equation by Benjamin~\cite{Benjamin1972Stability}. Its first application to the Klein-Gordon equation, resulting in orbital stability of standing waves, is due to Shatah~\cite{Shatah1983Stable}. Building on the work of Benjamin and others, the method was formalized by Grillakis, Shatah, and Straus, culminating in a comprehensive orbital stability theory~\cite{Grillakis1987Stability,Grillakis1990Stability} for traveling waves in general Hamiltonian systems with symmetries. A thorough introduction on the method and further references can be found in~\cite{kapitula_spectral_2013}.

Although the theory of Grillakis, Shatah, and Strauss provides an adequate framework for the stability analysis of solitary waves, its application to \emph{periodic} waves faces limitations. Specifically, the variational approach is heavily dependent on conserved quantities (such as the Hamiltonian), which require that the perturbed wave is either localized or periodic. As a result, the class of admissible perturbations for a periodic wave is restricted to those with the same period (co-periodic or harmonic perturbations) or those whose period is a multiple of that of the wave (subharmonic perturbations). In the Klein-Gordon equation~\eqref{e:KG}, orbital stability results of periodic waves against harmonic perturbations have been obtained for several classes of nonlinearities; see~\cite{Chen2023Orbital,Natali2008Stability,Palacios2022Orbital}. For an overview of co-periodic stability results in other Hamiltonian systems, as well as a general characterization of the Grillakis-Shatah-Strauss stability criteria in the periodic context, we refer to~\cite{Benzoni2013Stability,Bezoni2016Coperiodic}.

As outlined in~\cite{Benzoni2013Stability,Bezoni2016Coperiodic,Gallay2007Stability}, the nonlinear stability of periodic waves with respect to localized perturbations in Hamiltonian systems remains a largely unsolved problem, standing in sharp contrast with the dissipative case; see Remark~\ref{rem:dissipative}. To the authors' best knowledge, the only known nonlinear stability result to date for this class of perturbations concerns the orbital stability of plane waves in the one-dimensional nonlinear Schr\"odinger (NLS) equation 
\begin{align} \label{e:NLS}
\ri u_t + u_{xx} + f\big(|u|^2\big) u = 0, \qquad x, t \in \RR, \quad u(x,t) \in \mathbb{C}.
\end{align}
Thanks to the Galilean invariance of the NLS equation, the plane-wave solution~\eqref{e:plane-wave} can be transformed into the time-periodic, spatially homogeneous solution
\begin{align} \label{e:homosci}
u(x,t) = a\re^{\ri (\omega + k^2) t}
\end{align}
to~\eqref{e:NLS}. Its orbital stability with respect to localized perturbations was established by Zhidkov~\cite{Zhidkov2001Korteweg}. The result holds locally uniformly in space and is derived via variational methods: by taking a suitable formal linear combination of the classical conservation laws for mass and energy, one obtains a well-defined conserved quantity, which serves to close the nonlinear argument. As observed in~\cite{Gallay2007Stability}, this also implies orbital stability of the associated plane-wave solution under localized perturbations.

\begin{remark} \label{rem:dissipative}
The nonlinear stability theory of periodic waves in dissipative systems with respect to localized perturbations has advanced significantly over the last decades. The first nonlinear stability result was obtained for plane-wave solutions in the real Ginzburg-Landau equation~\cite{collet1992diffusive} by exploiting its gauge invariance. The use of mode filters in Bloch frequency domain~\cite{schneider1996Diffusive,Schneider1998Nonlinear}, combined with the modulational ansatz proposed in~\cite{doelman2009dynamics}, led to a breakthrough in the nonlinear stability analysis of periodic waves in large classes of dissipative systems such as reaction-diffusion models~\cite{Johnson2011Nonlinear,Johnson2013nonlocalized,sandstede2012diffusive} and systems of viscous conservation laws~\cite{Johnson2010nonlinear,johnson2014behavior}. Based on these developments, we believe that the orbital stability analysis of \emph{plane waves} in the complex Klein-Gordon equation presented in this paper offers a promising step toward developing an orbital stability theory for more general periodic waves in Hamiltonian systems.
\end{remark}

\subsection{Goal, challenges, and strategy} \label{sec:strategy}

Our aim is to develop a variational framework, in the spirit of~\cite{Zhidkov2001Korteweg}, for establishing orbital stability of plane waves in the complex Klein-Gordon equation~\eqref{e:KG} against localized perturbations. In contrast to the NLS equation where Galilean invariance enables a reduction to a spatially homogeneous solution to the NLS equation itself, no such reduction is available for the Klein-Gordon equation~\eqref{e:KG}.

To address this difficulty, we adopt an alternative reduction inspired by the spectral stability analysis~\cite{Demirkaya2015Spectral} of spatially periodic waves in~\eqref{e:KG}. Specifically, we introduce the new coordinate
\begin{align} \label{e:coord}
u(x,t) = w(x-ct,t)a\re^{\ri k x + \ri \omega t},
\end{align}
where $c \in \RR$ is a free parameter denoting the speed of the co-moving frame in which the dynamics of the new $w$-variable is observed. The Hamiltonian structure of~\eqref{e:KG} is preserved under the coordinate transformation~\eqref{e:coord}, and the plane-wave solution~\eqref{e:plane-wave} is mapped to the homogeneous equilibrium state $w(x,t) \equiv 1$. This transformation thus reduces the stability problem to that of a (nonlocalized) \emph{stationary} solution $\mathbf{w}_*$ to a Hamiltonian system of the form
\begin{align} \label{e:Hamsys2}
\mathbf{w}_t = \mathcal{J} \mathcal{H}'(\mathbf{w}).
\end{align}

The goal is to control the dynamics of a perturbed solution $\mathbf{w}(t)$ to~\eqref{e:Hamsys2} with initial condition $\mathbf{w}(0) = \mathbf{w}_* + \mathbf{z}_0$, where $\mathbf{z}_0$ is an $L^2$-localized perturbation. A crucial observation is that, since stationary solutions to~\eqref{e:Hamsys2} are critical points of the Hamiltonian, the formal expression 
\begin{align} \label{e:formalexpression}
E(\mathbf{z}(t)) = \mathcal{H}(\mathbf{w}_* + \mathbf{z}(t)) - \mathcal{H}(\mathbf{w}_*) = \mathcal{O}\big(\|\mathbf{z}(t)\|^2\big)
\end{align}
represents a well-defined conserved quantity for the $L^2$-localized perturbation $\mathbf{z}(t) = \mathbf{w}(t) - \mathbf{w}_*$. 

If the second variation $\mathcal{H}''(\mathbf{w}_*)$ were positive definite, then the energy acts as a Lyapunov functional and $\|\mathbf{z}(t)\|^2$ can be bounded in terms of $E(\mathbf{z}(t)) = E(\mathbf{z}_0)$, implying nonlinear stability. However, this is obstructed by the gauge symmetry of~\eqref{e:KG} which forces $0$ to lie in the spectrum of $\mathcal{H}''(\mathbf{w}_*)$. We emphasize that, since $\mathcal{H}''(\mathbf{w}_*)$ has constant coefficients, there is also no finite-codimensional subspace of $L^2(\RR)$ on which $\mathcal{H}''(\mathbf{w}_*)$ is positive definite, which is a key assumption in the classical Grillakis-Shatah-Strauss stability framework. 

To overcome this difficulty, we decompose the $w$-variable in polar coordinates 
\begin{align} \label{e:wpolar}
w(t) = \re^{\rho(t) + \ri \phi(t)},
\end{align}
where the phase variable $\phi(t)$ captures the neutral behavior caused by the gauge symmetry. This allows us to derive effective lower bounds on the conserved energy $E(\mathbf{z}(t))$, which control $\rho(t)$ and \emph{derivatives} of $\phi(t)$ over time. Our analysis shows that such bounds can be obtained provided that the spectral condition
\begin{align} \label{e:condition} a^2 f'\big(a^2\big) > 2\max\left\{0,-f\big(a^2\big)\right\}, \end{align} 
ensuring that $\mathcal{H}''(\mathbf{w}_*)$ is positive semi-definite, is satisfied. All in all, the obtained control on $\rho(t)$ and on derivatives of $\phi(t)$ is sufficient to close a nonlinear iteration argument, yielding an orbital stability result which is locally uniform in space.

\begin{remark}
The plane wave~\eqref{e:plane-wave} is a trivial example of a \emph{traveling-kink solution} 
\begin{align*} 
\psi(x,t) = \psi_0(x-st) \re^{\ri kx + \ri \omega t},
\end{align*}
where the wave speed $s$, wave number $k$, and temporal frequency $\omega$ are real parameters and the wave profile $\psi_0 \colon \RR \to \RR$ connects \emph{nonvanishing} asymptotic end states $\psi_\pm = \lim_{\xi \to \pm \infty} \psi(\xi) \neq 0$. When the wave profile is strictly monotone, i.e.~$\psi_0'(\xi) \neq 0$ for all $\xi \in \RR$, the use of variational methods to prove orbital stability is well-established in one-dimensional Hamiltonian systems. We refer to~\cite{Bethuel2008Orbital,Gerard2009Orbital,Henry1982Stability,Iliev1987Stability,Zhidkov1987cKG,Zhidkov1987NLS,Zhidkov2001Korteweg} for results in the (generalized) Korteweg-de Vries, Klein-Gordon, NLS, and nonlinear wave equations. A further extension to solitonic bubble solutions of the form $\psi(x,t) = \psi_0(x-st) \re^{\ri \phi_0(x-st)}$, with $\psi_0$ even and $(\psi_0(\xi),\phi_0'(\xi)) \to (\psi_\infty,0)$ as $\xi \to \pm \infty$ can be found in~\cite{Lin2002Stability}.
\end{remark}

\subsection{Main result}

We formulate our main result concerning the orbital stability of plane waves in the complex Klein-Gordon equation under localized perturbations. Since our analysis only requires control over the $L^2$-norm of the spatial and temporal derivatives of the phase variable, we are able to establish orbital stability under unbounded initial phase modulations; see Remark~\ref{rem:phase_offset}. Thus, we consider initial data of the form
\begin{align} \label{e:IC}
\begin{pmatrix}
u(x,0) \\ \partial_t u(x,0) 
\end{pmatrix} =
a\re^{\ri kx + \ri \theta_\infty(x)} \begin{pmatrix}
1 \\ \ri \omega
\end{pmatrix} + \begin{pmatrix}
w_0(x) \\ v_0(x)
\end{pmatrix}
\end{align}
where $\theta_\infty \colon \RR \to \RR$ is a continuously differentiable phase function with $\theta_\infty' \in H^1(\RR)$, and $w_0 \in H^1(\RR)$ and $v_0 \in L^2(\RR)$ are localized perturbations. 

We are now in a position to state the main result of this paper.

\begin{theorem}
\label{t:main}
Let $(a, k, \omega) \in (0,\infty) \times \RR^2$ with $(k,\omega) \neq (0,0)$. Assume that the existence condition~\eqref{e:disp-rel} and the spectral condition~\eqref{e:condition} are satisfied. Then, there exist constants $C,\delta > 0$ such that, whenever $\theta_\infty \in \cC^1(\RR,\RR)$ and $(w_0,v_0) \in H^2(\RR) \times H^1(\RR)$ satisfy $\theta_\infty' \in H^1(\RR)$ and
\begin{align}
\label{e:initial-smallness}
E_0 := \norm{\theta_\infty'}_{L^2(\RR)} + \norm{w_0}_{H^1(\RR)} + \norm{v_0}_{L^2(\RR)} \leq \delta,
\end{align}
there exists a unique global classical solution 
\begin{align} \label{e:regu}
u \in \cC\big([0,\infty),H^2_{\mathrm{loc}}(\RR)\big) \cap \cC^1\big([0,\infty),H^1_{\mathrm{loc}}(\RR)\big) \cap \cC^2\big([0,\infty),L^2_{\mathrm{loc}}(\RR)\big)
\end{align}
to the Klein-Gordon equation~\eqref{e:KG} with initial condition~\eqref{e:IC}. Furthermore:
\begin{itemize}
\item \emph{(Local orbital stability)} For every $R > 0$, $x_* \in \RR$, and $t \geq 0$ there exists a constant $\gamma = \gamma(x_*,t) \in \RR$ such that
\begin{align}\label{e:orbitalstab}
\norm{u(t) - a\re^{\ri k (\cdot) + \ri \omega t + \ri \gamma}}_{H^1(x_* - R,x_* + R)} \leq C\sqrt{R} E_0.
\end{align}
\item \emph{(Bounds on phase and radius)} There exist continuously differentiable functions $\varrho,\vartheta \in \cC^1\big(\RR \times [0,\infty),\RR\big)$ satisfying $\varrho(t), \vartheta(t) \in H^1(\RR)$ and $\varrho_t(t), \vartheta_t(t) \in L^2(\RR)$ such that the polar decomposition
\begin{align}
\label{polardecomp}
u(x,t) = a \re^{\ri kx + \ri \omega t + \ri \theta_\infty(x)} \, \re^{\varrho(x,t) + \ri \vartheta(x,t)}
\end{align}
holds for all $x \in \RR$ and $t \geq 0$. Moreover, we have the bounds
\begin{align}
\label{e:mainbounds}
\begin{split}
\norm{\varrho_t(t)}_{L^2(\RR)} + \norm{\vartheta_t(t)}_{L^2(\RR)} + \norm{\varrho(t)}_{H^1(\RR)} + \norm{\vartheta_x(t)}_{L^2(\RR)} & \leq C E_0,\\
\norm{\vartheta(t)}_{L^2(\RR)} &\leq CE_0 (t + 1)
\end{split}
\end{align}
for $t \geq 0$.
\end{itemize}
\end{theorem}

\begin{remark} \label{rem:phase_offset}
As an illustrative example for the initial phase modulation in Theorem~\ref{t:main}, consider
\begin{align*}
\theta_\infty(x) = k x_{-} \frac{\tanh(x) - 1}{2} - k x_{+} \frac{1 + \tanh(x)}{2}
\end{align*}
with $k,x_\pm \in \RR$. This function satisfies $\norm{\theta_\infty'}_{L^2(\RR)} = |k(x_+ - x_-)|/\sqrt{3}$. Hence, for any fixed $\delta>0$, it obeys $\|\theta_\infty'\|_{L^2(\RR)} < \delta$ provided that the phase offset $k(x_+ - x_-)$ is sufficiently small. The associated modulated plane wave $\smash{a\re^{\ri kx + \ri \theta_\infty(x)}}$ connects to the phase-shifted waves $\smash{a\re^{\ri k(x-x_\pm)}}$ at $\pm \infty$. Theorem~\ref{t:main} also allows for \emph{unbounded} initial phase modulations $\theta_\infty$. For instance, there exists a constant $C_0 > 0$ such that the unbounded smooth function $\theta_\infty(x) = \smash{(1 + \epsilon^4 x^2)^{\frac18}}$ satisfies $\|\theta_\infty'\|_{L^2(\RR)} \leq C_0\epsilon$ for all $\epsilon \in (0,1)$. Thus, for any fixed $\delta > 0$, the condition~\eqref{e:initial-smallness} can be met by taking $\epsilon \in (0,1)$ sufficiently small. 
\end{remark}

\subsection{Discussion and outlook}

Our spectral analysis in~\S\ref{sec:specstab} shows that Theorem~\ref{t:main} is sharp in the sense that it holds up to the spectral stability boundary. That is, at the critical value $a^2 f'(a^2) = 2 \max\{0,-f(a^2)\}$ where condition~\eqref{e:condition} ceases to be satisfied, the plane wave~\eqref{e:plane-wave} transitions from a spectrally stable to a spectrally unstable solution of~\eqref{e:KG}. Notably, the condition~\eqref{e:condition} only involves the leading-order terms in the expansion~\eqref{e:fexpans} of the nonlinearity. When $f(a^2) \geq 0$, the loss of stability corresponds to a transition from a defocusing to a focusing nonlinearity. In contrast, if $f(a^2) < 0$, the instability arises within the defocusing regime itself. We note that a change in stability is accompanied by a corresponding change in character (from elliptic to hyperbolic) of the Whitham's modulation equation, which can be derived via a long-wave ansatz. For further details and the rigorous justification of the Whitham's modulation equation for plane waves in the complex cubic Klein–Gordon equation, we refer to~\cite{hofbauer2022validity}.

Similar to the orbital stability result for time-periodic, spatially homogeneous solutions~\eqref{e:homosci} to the NLS equation in~\cite{Zhidkov2001Korteweg}, our result is locally uniform in space. Numerical simulations, see~\S\ref{sec:numerics}, suggest that an orbital stability result, where $H^1(x_*-R,x_*+R)$ in~\eqref{e:orbitalstab} is replaced by $H^1(\RR)$ or where the right-hand side of~\eqref{e:orbitalstab} is replaced by $CE_0$, is not attainable. Specifically, our simulations reveal that localized perturbations trigger a local phase rotation of the plane wave, which subsequently spreads with constant speed in both spatial directions. That is, perturbations initiate a two-sided expanding phase front. Our numerical analysis indicates that the $L^2$-norm of $\vartheta(t)$ grows algebraically at rate $\smash{t^{\frac12}}$, while its $L^\infty$-norm remains bounded. Thus, we do not expect the growth rate of the $L^2$-norm of $\vartheta(t)$ in Theorem~\ref{t:main} to be sharp. A compelling direction for future research would be to characterize the leading-order phase dynamics, which may, at least for long-wavelength perturbations, be governed by Whitham's modulation equation; see~\cite{hofbauer2022validity}. 

We note that the regularity assumptions on the initial data in Theorem~\ref{t:main} can be relaxed to $w_0 \in H^1(\RR,\CC)$ and $v_0 \in L^2(\RR,\CC)$ by working with mild rather than classical solutions of the perturbation equation; see equation~\eqref{e:short-time-dynamic} in~\S\ref{sec:local_existence}. One then obtains a function $u \in \cC\big([0,\infty),H^1_{\mathrm{loc}}(\RR)\big) \cap C^1([0,\infty),L^2_{\mathrm{loc}}(\RR)\big)$ which satisfies~\eqref{e:orbitalstab},~\eqref{polardecomp} and~\eqref{e:mainbounds}, and is the limit of a sequence of classical solutions $\{u_n(t)\}_n$ to~\eqref{e:KG} with initial data
\begin{align*}
\begin{pmatrix}
u_n(x,0) \\ \partial_t u_n(x,0) 
\end{pmatrix} =
a\re^{\ri kx + \ri \theta_\infty(x)} \begin{pmatrix}
1 \\ \ri \omega
\end{pmatrix} + \begin{pmatrix}
w_n(x) \\ v_n(x)
\end{pmatrix},
\end{align*}
where $\{(w_n,v_n)^\top\}_n$ is a sequence in $H^2(\RR) \times H^1(\RR)$ converging to $(w_0,v_0)^\top$ in $H^1(\RR) \times L^2(\RR)$. For the sake of simplicity of exposition, we do not pursue this generalization here.

We expect that the variational approach in this paper can be extended to higher spatial dimensions. In such settings, an intriguing direction is to explore whether the method can be combined with dispersive estimates to derive optimal or sharper bounds on the $\varrho$- and $\vartheta$-variables in Theorem~\ref{t:main}. Owing to stronger dispersive decay in higher dimensions, one could even try to establish an \emph{asymptotic} nonlinear stability result for plane-wave solutions. To the best of the authors' knowledge, no asymptotic nonlinear stability results currently exist for periodic waves in dispersive Hamiltonian systems; we refer to~\cite{Audiard2024Linear,Rodrigues2018Linear} for \emph{linear} asymptotic stability results for periodic solutions to the (generalized) Korteweg-de Vries equation. By contrast, the asymptotic nonlinear stability theory for soliton and kink solutions in such systems has seen significant advances over the past decades; see, for instance, the surveys~\cite{Cuccagna2021Survey,Kowalczyk2017Asymptotic} and references therein. Finally, another interesting avenue for future research is whether the orbital stability analysis presented here can be extended beyond plane waves to other classes of periodic solutions in Hamiltonian systems, particularly those that do not permit a reduction to a spatially homogeneous state.

\subsection{Organization} In~\S\ref{sec:coord_trans}, we apply the coordinate transform~\eqref{e:coord} to the Klein-Gordon equation and use the Hamiltonian structure of the transformed system to establish an energy. In~\S\ref{sec:local_existence}, we introduce the perturbation of the modulated plane wave, derive the corresponding perturbation equation, and establish local existence and uniqueness of its solutions. In~\S\ref{sec:energy}, we show that the energy is conserved and derive associated lower and upper bounds. Section~\ref{sec:mainresultproof} is devoted to the proof of our main result, Theorem~\ref{t:main}. Numerical simulations, supporting our analysis, are presented in~\S\ref{sec:numerics}. Finally, we carry out a spectral analysis in~\S\ref{sec:specstab}, demonstrating that our orbital stability result holds up to the spectral stability boundary.

\subsection{Data availability statement}

The numerical simulations displayed in Figures \ref{fig1}, \ref{fig2} and \ref{fig3} have been obtained using Matlab (Version 24.1.0, R2024a), and the code is available through the link \url{https://www.waves.kit.edu/downloads/CRC1173_Preprint_2025-23_Codes.zip}.

\subsection{Acknowledgments} This work is funded by the Deutsche Forschungsgemeinschaft (DFG, German Research Foundation) -- Project-ID 258734477 -- SFB 1173. We would like to thank Maximilian Ruff for helping us set up the numerical simulations.

\section{Coordinate transformation and energy} \label{sec:coord_trans}

In this section, we apply the coordinate transformation~\eqref{e:coord}  to the complex Klein–Gordon equation~\eqref{e:KG}, under which the plane-wave solution~\eqref{e:plane-wave} is mapped to the equilibrium state $w(x,t) \equiv 1$. We establish that the transformed system retains a Hamiltonian structure, from which we derive an energy using the strategy outlined in~\S\ref{sec:coord_trans}. The conservation of this energy forms the backbone of our orbital stability argument.

Inserting~\eqref{e:coord} into~\eqref{e:KG} and using~\eqref{e:disp-rel}, we find that the evolution of the $w$-variable is governed by the equation
\begin{align} \label{e:KGw}
\begin{split}
w_t &= v,\\
v_t &= \left(1 - c^2\right) w_{yy} + 2c v_y - 2 \ri \omega v + 2 \ri\left(c\omega + k\right) w_y - \left(f\big(a^2|w|^2\big) - f\big(a^2\big)\right) w.
\end{split}
\end{align}
In $(w,v)$-coordinates the plane-wave solution~\eqref{e:plane-wave} corresponds to the homogeneous rest state $(1,0)^\top$. In terms of the real coordinate vector $\mathbf{w} = (\Real(w),\Imag(w),\Real(v),\Imag(v))^\top$, the transformed system~\eqref{e:KGw} takes the Hamiltonian form~\eqref{e:Hamsys2}, where $\mathcal{J}$ denotes the skew-symmetric operator 
\begin{align*}
\mathcal{J} = \begin{pmatrix} 0 & 0 & 1 & 0 \\ 0 & 0 & 0 & 1 \\ -1 & 0 & 2 c \partial_y & 2 \omega \\ 0 & -1 & -2\omega & 2c\partial_y \end{pmatrix},
\end{align*}
$\mathcal{H} \colon \mathcal{X} \to \RR$ is the Hamiltonian
\begin{align*}
\mathcal{H}(\xi,\zeta,\mu,\nu) &= \frac{1}{2} \int_{\mathcal{\RR}} (1-c^2) \left(\xi'(y)^2 + \zeta'(y)^2\right) + \mu(y)^2 + \nu(y)^2 + \frac{1}{a^2} F\left(a^2 \big(\xi(y)^2 + \zeta(y)^2\big)\right)  \\ 
&\qquad \qquad - f\big(a^2\big) \left(\xi(y)^2 + \zeta(y)^2\right) + 2(c\omega + k)\left(\xi(y)\zeta'(y) - \xi'(y) \zeta(y)\right)   \de y,
\end{align*}
and $\mathcal{X} = H^1(\RR) \times H^1(\RR) \times L^2(\RR) \times L^2(\RR)$ is the associated Hilbert space. A nonlinear functional, corresponding to the formal expression~\eqref{e:formalexpression}, is then given by the energy
\begin{align} \label{e:defenergy2}
\begin{split}
E_c(w,v) &= \frac{1}{2} \int_{\mathcal{\RR}} (1-c^2) |w'(y)|^2 + |v(y)|^2 + U\big(|w(y)|^2\big) - 2(c\omega + k) \Imag\left(w(y) \overline{w'(y)}\right)\de y
\end{split}
\end{align}
where we denote
\begin{align} \label{defU}
U(s) = \int_1^s f\big(a^2 \upsilon\big) - f\big(a^2\big) \, \de \upsilon.
\end{align}
For the orbital stability analysis of the equilibrium state $(1,0)^\top$ in~\eqref{e:KGw}, it is convenient to use the polar representation
\begin{align*}
\begin{pmatrix} w(y,t)\\ v(y,t)\end{pmatrix} = \begin{pmatrix}  \re^{\rho(y,t) + \ri \phi(y,t)}\\ v(y,t)\end{pmatrix}.
\end{align*}
Then, $\rho$ and $v$ measure the deviation from the manifold of equilibria $\mathcal{M} = \{(\re^{\ri \gamma},0)^\top : \gamma \in \RR\}$ in~\eqref{e:KGw}, and $\phi$ tracks the solution's orbit along $\mathcal{M}$ induced by gauge symmetry. Expressing
\begin{align} \label{e:defenergy}
\begin{split}
E_c(w,v)  &= \frac{1}{2} \int_{\mathcal{\RR}} (1-c^2) \left(\rho'(y)^2 + \phi'(y)^2\right) \re^{2\rho(y)} + |v(y)|^2 + U\big(\re^{2\rho(y)}\big)\\
&\qquad \qquad + 2(c\omega + k) \left(\phi'(y) \left(\re^{2\rho(y)} - 1\right) + \phi'(y)\right) \de y,
\end{split}
\end{align}
it follows that the energy $E_c(w,v)$ is well-defined as long as $v \in L^2(\RR)$ and $w = \re^{\rho + \ri \phi}$ with $\rho,\phi \in H^1(\RR)$. In~\S\ref{sec:energy}, we slightly modify this energy to allow for unbounded initial phase modulations. Moreover, we prove that this modified energy is conserved and use~\eqref{e:condition} and~\eqref{e:defenergy} to derive lower and upper bounds which are sufficient to close an orbital stability argument.

\section{Local existence analysis} \label{sec:local_existence}

In this section, we introduce a localized perturbation of the modulated plane wave, derive an associated perturbation equation, and establish local existence and uniqueness of its solutions. We then transfer the obtained localization and regularity properties to the polar representation of the perturbation, which, as outlined in~\S\ref{sec:strategy}, effectively captures the neutral behavior induced by gauge symmetry in our orbital stability analysis. 

The coordinate transformation~\eqref{e:coord} maps the modulated plane wave
\begin{align*}
u(x,t) = a \re^{\ri k x + \ri \omega t + \ri \theta_\infty(x - c t)}
\end{align*}
to $w(y,t) = \re^{\ri \theta_\infty(y)}$. We measure the deviation from the modulated plane wave by writing
\begin{align} \label{e:def_zvariable}
w(y,t) = \re^{\ri \theta_\infty(y)} \left(1 + z(y,t)\right). 
\end{align}
Our next step is to derive an evolution equation for the perturbation $z$ and its temporal derivative by inserting~\eqref{e:def_zvariable} into~\eqref{e:KGw}. It is convenient to rescale the temporal derivative by $1/\sqrt{1-c^2}$, which renders the principal linear part of the equation skew-adjoint. Accordingly, we assume that the free variable $c \in \RR$ satisfies $c \in (-1,1)$, and introduce
\begin{align}
Z(y,t) = \begin{pmatrix} z(y,t) \\ \frac{z_t(y,t)}{\sqrt{1-c^2}}\end{pmatrix}. \label{e:defZ}
\end{align}
We arrive at the semilinear evolution system
\begin{align} \label{e:short-time-dynamic}
Z_t = \Lambda Z + N(Z)
\end{align}
posed on the Hilbert space $X = H^1(\RR) \times L^2(\RR)$, where the linear operator $\Lambda \colon D(\Lambda) \subset X \to X$ is defined on the domain $D(\Lambda) = H^2(\RR) \times H^1(\RR)$ by
\begin{align*}
\Lambda = \sqrt{1 - c^2}
\begin{pmatrix}
0 & 1 \\
\partial_{yy} - 1 & \frac{2c}{\sqrt{1-c^2}} \partial_y
\end{pmatrix},
\end{align*}
and the nonlinearity $N \colon X \to X$ is given by
\begin{align*}
N(Z) = \begin{pmatrix} 0 \\ \widetilde{N}(Z)\end{pmatrix}
\end{align*}
with
\begin{align*}
\widetilde{N}(z_1,z_2) &= \frac{2\ri(k + c \omega)}{\sqrt{1-c^2}} \partial_y z_1 - 2\ri \omega z_2 + \sqrt{1-c^2} \, z_1 - \frac{1 + z_1}{\sqrt{1-c^2}} \bigg(f\big(a^2\absolute{1 + z_1}^2\big) - f\big(a^2\big)\bigg)\\
&\qquad + 2 \ri \theta_\infty' \Big(\sqrt{1-c^2} \, \partial_y z_1 + c z_2\Big) + (1 + z_1) \bigg(\sqrt{1 - c^2} \big(\ri \theta_\infty'' - (\theta_\infty')^2\bigg) - \frac{2 (k + c \omega)}{\sqrt{1-c^2}} \theta_\infty'\bigg).
\end{align*}

Local existence and uniqueness of classical solutions to~\eqref{e:short-time-dynamic} follows from standard semigroup theory.

\begin{proposition} \label{p:short-time-existence}
Assume that $c \in (-1,1)$. Let $\theta_\infty \in \cC^1(\RR)$ be such that $\theta_\infty' \in H^1(\RR)$. Take $Z_0 \in D(\Lambda)$. Then, there exist a maximal time $T = T(Z_0) \in (0,\infty]$ and a unique classical solution $Z \in \cC([0,T),X) \cap \cC^1([0,T),D(\Lambda))$ to~\eqref{e:short-time-dynamic} with initial value $Z(0) = Z_0$. Moreover, if $T < \infty$, then it holds
\begin{align} 
\label{e:blow-up}
\limsup_{t \uparrow T} \, \norm{Z(t)}_X = \infty.
\end{align}
\end{proposition}
\begin{proof}
Note that $\Lambda$ is a skew-adjoint operator. As such, it generates a unitary group on the Hilbert space $X$ by Stone's Theorem; see~\cite[Corollary II.3.24]{engel_one-parameter_2000}. On the other hand, using that $H^1(\RR)$ continuously embeds into $L^\infty(\RR)$ and that we have $f \in \cC^1(\RR)$, one readily infers that $N$ is a well-defined locally Lipschitz continuous map. The result directly follows from standard semigroup theory for semilinear evolution equations; see~\cite[Theorem~4.3.4 and~Proposition~4.3.9]{cazenave1998introduction}. 
\end{proof}

In our orbital stability analysis, we employ the polar representation 
\begin{align}
\label{e:polar-decomposition} 
z(y,t) = \eu^{\rho(y,t)+\iu \theta (y,t)} - 1
\end{align}
of the perturbation to control its long-term dynamics. The following result shows that $\rho$ and $\theta$ inherit the localization and regularity properties of $z$, as long as $1+z(t)$ remains sufficiently close to the circle of radius $1$ centered at the origin. 

\begin{proposition}
\label{p:polar-coordinate}
Take $c \in (-1,1)$, $\delta_1 \in (0,\frac14)$, and $\theta_\infty \in \cC^1(\RR,\RR)$ with $\theta_\infty' \in H^1(\RR)$. Let $\rho_0, \theta_0 \in H^2(\RR,\RR)$, $z_1 \in H^2(\RR,\CC)$ and $z_2 \in H^1(\RR,\CC)$ be such that
\begin{align} \label{e:z0_annulus}
\big\lVert \re^{\rho_0} - 1\big\rVert_{L^\infty(\RR)} \leq \delta_1, \qquad z_1 = \re^{\rho_0 + \ri \theta_0}-1.
\end{align} 
Denote by $z(t)$ the first coordinate of the classical solution $Z \colon [0,T(Z_0)) \to D(\Lambda)$ to~\eqref{e:short-time-dynamic} with initial condition $Z(0) = Z_0 = (z_1,z_2) \in D(\Lambda)$, which was established in Proposition~\ref{p:short-time-existence}. 

Then,
\begin{align*}
t_{\max} = \sup \Set{t \in [0, T(Z_0)) : \big\lVert\, \absolute{1 + z(t)} - 1\big\rVert_{L^\infty(\RR)} \leq 2 \delta_1}
\end{align*}
is strictly positive, and there exist unique continuously differentiable functions $\rho,\theta \in \cC^1\big(\RR \times [0,t_{\max}),\RR\big)$ with $\rho(0) = \rho_0$ and $\theta(0) = \theta_0$ such that the polar decomposition~\eqref{e:polar-decomposition} holds for all $y \in \RR$ and $t\in [0, t_{\max})$. Moreover, we have $\rho(t),\theta(t) \in H^1(\RR)$ and $\rho_t(t),\theta_t(t) \in L^2(\RR)$ with
\begin{align}\label{e:bound-rho}
\begin{split}
\norm{\rho(t)}_{L^\infty(\RR)} \leq 4 \delta_1, \qquad \norm{\rho(t)}_{L^2(\RR)} \leq 4\norm{z(t)}_{L^2(\RR)},
\end{split}
\end{align}
and
\begin{align} \label{e:bound-theta}
\begin{split}
\norm{\rho_y(t)}_{L^2(\RR)}^2 + \norm{\theta_y(t)}_{L^2(\RR)}^2 &\leq \frac{\norm{z_y(t)}_{L^2(\RR)}^2}{(1 - 2 \delta_1)^2},\\
\norm{\rho_t(t)}_{L^2(\RR)}^2 + \norm{\theta_t(t)}_{L^2(\RR)}^2 &\leq \frac{\norm{z_t(t)}_{L^2(\RR)}^2}{(1 - 2 \delta_1)^2},\\
\norm{\theta(t)}_{L^2(\RR)} &\leq \norm{\theta_0}_{L^2(\RR)} + t \frac{\sup_{s \in [0,t]} \norm{\partial_s z(s)}_{L^2(\RR)}}{1-2\delta_1}
\end{split}
\end{align}
for all $t\in [0, t_{\max})$.
\end{proposition}
\begin{proof}
First, we recall from~\eqref{e:defZ} that the second coordinate of $Z(t)$ is given by the temporal derivative $z_t(t)/\sqrt{1-c^2}$. Therefore, Proposition~\ref{p:short-time-existence} ensures that
\begin{align} \label{e:regz}
z \in \cC\big([0,T),H^2(\RR)\big) \cap \cC^1\big([0,T),H^1(\RR)\big) \cap \cC^2\big([0,T),L^2(\RR)\big).
\end{align}
In particular, using the continuous embedding $H^1(\RR) \hookrightarrow L^\infty(\RR)$, this implies that $z \colon \RR \times [0,T) \to \CC$ is continuously differentiable in both space and time. Moreover, in conjunction with the condition~\eqref{e:z0_annulus}, it guarantees that $t_{\max}$ is strictly positive. 

Expanding the real logarithm at $1$ as  $\log(x)=x-1+\mathcal{O}((x-1)^2)$, we find that, for any $\psi \in \CC$ with $\Absolute{\absolute{1 + \psi} - 1} \leq 2\delta_1$, we have 
\begin{align*}
\log\absolute{1 + \psi} \leq 2 \absolute{\absolute{1 + \psi} - 1}.
\end{align*}
Hence, the function $\rho \colon \RR \times [0,t_{\max}) \to \RR$ given by $\rho(y,t) = \log \absolute{1 + z(y,t)}$ is well-defined, continuously differentiable in space and time, and satisfies $\rho(\cdot,0) = \rho_0$ and $\smash{\norm{\rho(t)}_{L^\infty(\RR)}  \leq 4 \delta_1} < 1$ for $t \in [0,t_{\max})$, which yields the first bound in~\eqref{e:bound-rho}.  Furthermore, the definition of $t_{\max}$ and the equality $\smash{\re^{\rho(t)} = \absolute{1 + z(t)}}$ ensure that for all $t\in [0, t_{\max})$ we have
\begin{align}
\label{e:bound-exp-rho}
1 - 2\delta_1 \leq \inf_{y \in \RR} \re^{\rho(y,t)}.
\end{align}
On the other hand, since the function $\varsigma \colon \RR \times [0, t_{\max}) \to \SS^1$ given by $\varsigma(y,t) = \re^{-\rho(y,t)}(1+z(y,t))$ is continuously differentiable in space and time, there exists a unique argument function $\theta \in \cC^1\big(\RR \times [0,t_{\max}), \RR\big)$ with $\theta(\cdot,0) = \theta_0$, which satisfies~\eqref{e:polar-decomposition} for all $t \in [0,t_{\max})$.

We now prove that $\rho$, $\theta$, and their spatial and temporal derivatives are localized in space. To this end, we note that $\absolute{\eu^{r} - 1} \geq \frac14 \absolute{r}$ holds for all $r \in (-1,1)$ by the mean value theorem. Therefore, using~\eqref{e:bound-exp-rho} and the fact that $\|\rho(t)\|_{L^\infty(\RR)} \leq 1$, we obtain
\begin{align*}
\begin{split}
\norm{z(t)}_{L^2(\RR)} &= \norm{\eu^{\rho(t)+\iu \theta(t)}-1}_{L^2(\RR)} \geq  \norm{\eu^{\rho(t)}-1}_{L^2(\RR)} \geq \frac{1}{4} \norm{\rho(t)}_{L^2(\RR)}, \\
\norm{z_y(t)}_{L^2(\RR)}^2 &= \norm{(\rho_y(t)+\iu\theta_y(t))\eu^{\rho(t)+\iu \theta(t)}}_{L^2(\RR)}^2 = \norm{\eu^{\rho(t)} {\rho_y(t)}}_{L^2(\RR)}^2 + \norm{\re^{\rho(t)} {\theta_y(t)}}_{L^2(\RR)}^2 \\
&\geq  (1 - 2\delta_1)^2 \left(\norm{{\rho_y(t)}}_{L^2(\RR)}^2  + \norm{{\theta_y(t)}}_{L^2(\RR)}^2 \right).
\end{split}
\end{align*}
for $t \in [0,t_{\max})$, implying the second bound in~\eqref{e:bound-rho} and the first bound in~\eqref{e:bound-theta}. Thus, we find that $z(t) \in H^1(\RR)$ implies $\rho(t) \in H^1(\RR)$ and $\theta_y(t) \in L^2(\RR)$ for all $t \in [0,t_{\max})$. Similarly,~\eqref{e:bound-exp-rho} leads to 
\begin{align*}
\begin{split}
\norm{z_t(t)}_{L^2(\RR)}^2 &= \norm{(\rho_t(t)+\iu\theta_t(t))\eu^{\rho(t)+\iu \theta(t)}}_{L^2(\RR)}^2\\ 
&\geq (1 - 2 \delta_1)^2 \left(\norm{\rho_t(t)}_{L^2(\RR)}^2 + \norm{\theta_t(t)}_{L^2(\RR)}^2\right),
\end{split}
\end{align*}
for $t \in [0,t_{\max})$, yielding the second bound in~\eqref{e:bound-theta}. Finally, the above estimate, together with the fundamental theorem of calculus
\begin{align*} 
\theta(y,t) = \theta(y,0) + \int_0^t \partial_s \theta(y,s) \dd s, \qquad y \in \RR, \, t \in [0,t_{\max})
\end{align*}
provide the last claimed bound in~\eqref{e:bound-theta}.
\end{proof}

\section{Energy estimates} \label{sec:energy}

In this section, we slightly modify the energy~\eqref{e:defenergy2} and establish upper and lower bounds, which play a pivotal role in our stability argument. Expressing~\eqref{e:defenergy2} in $Z$-coordinates via~\eqref{e:def_zvariable} and~\eqref{e:defZ} and correcting the integrand with $-2(k+c\omega) \theta_\infty'$, we arrive at the quantity
\begin{align}
\label{e:energy}
\begin{split}
E(z_1, z_2) &= \frac{1}{2} \int_\RR (1 - c^2) \left(\absolute{z_2(y)}^2 + \Absolute{\partial_y \left((1 + z_1(y)) \re^{\ri \theta_\infty(y)}\right)}^2\right) + U\left(\absolute{1 + z_1(y)}^2\right) \\
&\qquad - \, 2(k + c\omega) \left(\Imag \left((1 + z_1(y)) \overline{z_1'(y)}\right) - \theta_\infty'(y) \absolute{1 + z_1(y)}^2 + \theta_\infty'(y) \right) \dd y,
\end{split}
\end{align}
where we recall that $U \in \cC^3\big([0,\infty),\RR\big)$ is given by~\eqref{defU}. The correction term is chosen so as to eliminate linear terms in $\theta_\infty'$ from the integrand, which need not be integrable under the sole assumption $\theta_\infty' \in L^2(\RR)$. Before showing that the energy is conserved and bounded from below, we first prove that it is indeed well-defined and bounded from above.

\begin{lemma} \label{l:finite-energy}
Let $c \in (-1,1)$. There exists $C_1 > 0$ such that for each $\rho,\theta \in H^1(\RR,\RR)$, $z_1 \in H^1(\RR,\CC)$, $z_2 \in L^2(\RR,\CC)$, and $\theta_\infty \in \cC^1(\RR,\RR)$ satisfying $\theta_\infty' \in L^2(\RR)$ and
\begin{align*}
z_1 = \re^{\rho + \ri \theta} - 1, \qquad \norm{\rho}_{L^\infty(\RR)} \leq 1,
\end{align*} 
the energy $E(z_1,z_2)$ is well-defined and enjoys the upper bound
\begin{align*}
\absolute{E(z_1,z_2)} \leq C_1 \left(\norm{z_2}_{L^2(\RR)}^2 + \norm{\rho}_{H^1(\RR)}^2 + \norm{\theta'}_{L^2(\RR)}^2 + \norm{\theta_\infty'}_{L^2(\RR)}^2\right).
\end{align*}
\end{lemma}
\begin{proof}
We successively bound the three terms in the expression for $E$. First, using Young's inequality, we establish 
\begin{align*}
\Absolute{\partial_y \left((1 + z_1(y)) \re^{\ri \theta_\infty(y)}\right)}^2 &\leq 2\absolute{z_1'(y)}^2 + 2\absolute{1 + z_1(y)}^2 \absolute{\theta_\infty'(y)}^2\\ 
&= 2\left(\rho'(y)^2 + \theta'(y)^2 + \theta_\infty'(y)^2\right) \re^{2\rho(y)}
\end{align*}
for $y \in \RR$. Combining this pointwise bound with $\norm{\re^{2\rho}}_{L^\infty(\RR)} \leq \re^{2}$ yields the desired estimate on the first term in~\eqref{e:energy}. 

For the second term, we note that $U(1) = U'(1) = 0$. Hence, Taylor's Theorem yields a constant $C_* > 0$ such that 
\begin{align} \label{e:taylorMVT}
|\re^{2s} - 1| \leq C_* |s|, \qquad |U(\re^{2s})| \leq C_*|s|^2,
\end{align}
for $s \in [-1,1]$. As a consequence, it holds
\begin{align*}
U\left(\absolute{1 + z_1(y)}^2\right) = U\big(\re^{2\rho(y)}\big) \leq C_*|\rho(y)|^2
\end{align*}
for $y \in \RR$, which establishes the desired bound on the second term in~\eqref{e:energy}. 

For the third term, we compute
\begin{align} \label{e:computation}
\begin{split}
&\Imag\left((1 + z_1(y)) \overline{z_1'(y)}\right) - \theta_\infty'(y) \absolute{1 + z_1(y)}^2 + \theta_\infty'(y)\\ 
&\qquad = \theta_\infty'(y) -\big(\theta_\infty'(y)+\theta'(y)\big)\re^{2\rho(y)} = -\theta'(y) - \left(\theta_\infty'(y) + \theta'(y)\right) \left(\re^{2\rho(y)} - 1\right)
\end{split}
\end{align}
for $y \in \RR$. Hence, using that $\theta \in H^1(\RR)$ and applying the estimate~\eqref{e:taylorMVT}, we infer
\begin{align*}
&|k+c\omega|\Absolute{\int_\RR \Imag\left((1 + z_1(y)) \overline{z_1'(y)}\right) - \theta_\infty'(y) \absolute{1 + z_1(y)}^2 + \theta_\infty'(y) \dd y}\\ &\qquad = |k+c\omega|\Absolute{\int_\RR \left(\theta_\infty'(y) + \theta'(y)\right) \left(\re^{2\rho(y)} - 1\right) \dd y} \\
&\qquad \leq C_* \left( \norm{\theta'}_{L^2(\RR)} + \norm{\theta_\infty'}_{L^2(\RR)}\right)\norm{\rho}_{L^2(\RR)},
\end{align*}
which yields the desired bound on the third term in~\eqref{e:energy} by evoking Young's inequality. This completes the proof of the lemma.
\end{proof}

Next, we show that the energy is conserved under the evolution of~\eqref{e:short-time-dynamic}.

\begin{lemma} \label{l:conservation-energy}
Let $c \in (-1,1)$ and $\theta_\infty \in \cC^1(\RR,\RR)$ be such that $\theta_\infty' \in H^1(\RR)$. Take $Z_0 \in D(\Lambda)$. Let $Z \in \cC([0,T),X) \cap \cC^1([0,T),D(\Lambda))$ be the classical solution to~\eqref{e:short-time-dynamic} with initial condition $Z(0) = Z_0$, which was established in Proposition~\ref{p:short-time-existence}. Let $\tau \in (0,T]$ be such that the energy $E(Z(t))$ is well-defined for all $t \in [0,\tau)$. Then, we have
\begin{align*}
E(Z(t)) = E(Z_0), \qquad t \in [0,\tau).
\end{align*}
\end{lemma}
\begin{proof}
Denote $Z(t) = \transp{(z(t),v(t))}$. We successively differentiate the different terms in the expression~\eqref{e:energy} for the energy with respect to time. First, using $z_t(t) = \sqrt{1-c^2} \, v(t)$, we compute
\begin{align*}
\frac{\dd }{\dd t} \frac{1}{2} \int_\RR (1 - c^2) \Absolute{v(y,t)}^2 \dd y = \sqrt{1-c^2} \Real{\Scalp{v_t(t), z_t(t)}_{L^2(\RR)}}
\end{align*}
for $t \in [0,T)$. Moreover, using integration by parts, we derive
\begin{align*}
\frac{\dd }{\dd t} \left( \frac{1}{2} \int_\RR \Absolute{\partial_y \left((1 + z(y,t)) \re^{\ri \theta_\infty(y)}\right)}^2 \dd y \right) &=\Real{\Scalp{\partial_y \left((1 + z(t) \re^{\ri \theta_\infty}\right), \partial_y \left(z_t(t) \re^{\ri \theta_\infty}\right)}_{L^2(\RR)}}\\
&=-\Real{\Scalp{\re^{-\ri \theta_\infty} \partial_{yy} \left((1 + z(t)) \re^{\ri \theta_\infty}\right), z_t(t)}_{L^2(\RR)}}.
\end{align*}
for $t \in [0,T)$. We proceed with differentiating the second term in the energy. Recalling the definition~\eqref{defU} of $U$, we compute
\begin{align*}
\frac{\dd}{\dd t} \left(\frac{1}{2} \int_\RR  U\left(\absolute{1 + z(y,t)}^2\right) \dd y\right) = \Real{\Scalp{(1 + z(t))\Big(f\left(a^2 \absolute{1 + z}^2\right) - f\left(a^2\right)\Big), z_t(t)}_{L^2(\RR)}}
\end{align*}
for $t \in [0,T)$. For the third term, we use integration by parts to arrive at 
\begin{align*}
\frac{\dd}{\dd t} \frac{1}{2} \int_\RR \Imag\left((1 + z(y,t)) \overline{z_y(y,t)}\right) \dd y &= \frac12 \Imag\Big(\Scalp{z_t(t),z_y(t)}_{L^2(\RR)} + \Scalp{1+z(t),z_{yt}(t)}_{L^2(\RR)}\Big)\\ 
&= \frac12 \Imag\Big(\Scalp{z_t(t),z_y(t)}_{L^2(\RR)} - \Scalp{z_y(t),z_t(t)}_{L^2(\RR)}\Big)\\ 
&= \Real\Big(\ri \Scalp{z_y(t),z_t(t)}_{L^2(\RR)}\Big)
\end{align*}
for $t \in [0,T)$. Moreover, we find
\begin{align*}
\frac{\dd}{\dd t} \frac{1}{2} \int_\RR\theta_\infty'(y) \absolute{1 + z(y,t)}^2 - \theta_\infty'(y)\dd y = \Real{\Scalp{\theta_\infty' (1 + z(t)), z_t(t)}_{L^2(\RR)}}
\end{align*}
for $t \in [0,T)$. Gathering the above computations, using that $Z(t) = \transp{(z(t),v(t))}$ is a classical solution of~\eqref{e:short-time-dynamic}, recalling $z_t(t) = \sqrt{1-c^2} \, v(t)$, and integrating by parts, we obtain
\begin{align*}
\frac{\dd}{\dd t} E(Z(t)) &= \Real{\Scalp{2c z_{ty}(t) - 2 \ri (\omega - c\theta_\infty') z_t(t), z_t(t)}_{L^2(\RR)}}\\
&= c \Real{\Scalp{z_{ty}(t), z_t(t)}_{L^2(\RR)}} - c \Real{\Scalp{z_{t}(t), z_{ty}(t)}_{L^2(\RR)}} = 0
\end{align*}
for $t \in [0,T)$, which finishes the proof.
\end{proof}

We conclude this section by proving that the energy is bounded from below, provided that the spectral condition~\eqref{e:condition} holds.

\begin{lemma} \label{l:energy-controls-norm}
Let $(a,k,\omega) \in (0,\infty) \times \RR^2$ with $(k,\omega) \neq (0,0)$. Assume that the existence condition~\eqref{e:disp-rel} and the spectral condition~\eqref{e:condition} hold. Then, there exist constants $C_2 > 0$, $\delta_2 \in (0,1)$, and $c \in (-1,1)$ such that for each $\rho,\theta \in H^1(\RR,\RR)$, $z_1 \in H^1(\RR,\CC)$, $z_2 \in L^2(\RR,\CC)$, and $\theta_\infty \in \cC^1(\RR,\RR)$ such that $\theta_\infty' \in L^2(\RR)$ and
\begin{align*}
z_1 = \re^{\rho + \ri \theta} - 1, \qquad \|\rho\|_{L^\infty(\RR)} \leq \delta_2,
\end{align*}
the energy $E(z_1,z_2)$ enjoys the lower bound
\begin{align}
\label{e:energy-bound-below}
\norm{z_2}_{L^2(\RR)}^2 + \norm{\rho}_{H^1(\RR)}^2 + \norm{\theta'}_{L^2(\RR)}^2 \leq C_2 \left(E(z_1,z_2) + \norm{\theta_\infty'}_{L^2(\RR)}^2\right),
\end{align}
where we note that $k + c\omega = 0$ whenever $f(a^2) > 0$.
\end{lemma}
\begin{proof}
We successively bound the terms in the expression~\eqref{e:energy} for the energy from below. First, we note that $U(1) = U'(1) = 0$ and $U''(1) = a^2 f'(a^2)$. Hence, Taylor's Theorem affords a constant $C_* > 0$ such that
\begin{align} \label{e:Taylor2}
\left|U\big(\re^{2s}\big) - 2 a^2 f'(a^2) s^2\right| \leq C_* |s|^3, \qquad
\left|\re^{2s} - 1 - 2s\right| \leq C_*|s|^2, \qquad \left|\re^{2s} - 1\right| \leq C_* |s|
\end{align}
for $s \in [-1,1]$. Taking $\delta_2 \in (0,1)$ such that $1-C_*\delta_2 > 0$, we use~\eqref{e:Taylor2} to establish the pointwise lower bound
\begin{align*}
\Absolute{\partial_y\left( (1 + z_1(y)) \re^{\ri \theta_\infty(y)}\right)}^2 &= \re^{2\rho(y)} \left(\absolute{\rho'(y)}^2 + \absolute{\theta'(y) + \theta_\infty'(y)}^2\right)\\ 
&\geq (1 - C_*\delta_2) \left(\absolute{\rho'(y)}^2 + \absolute{\theta'(y) + \theta_\infty'(y)}^2\right)
\end{align*}
for $y \in \RR$. Hence, we arrive at 
\begin{align} \label{e:first_lower}
\begin{split}
&\int_\RR \Absolute{\partial_y\left( (1 + z_1(y)) \re^{\ri \theta_\infty(y)}\right)}^2 \dd y \geq (1 - C_* \delta_2) \left(\norm{\rho'}_{L^2(\RR)}^2 + \norm{\theta' + \theta_\infty'}_{L^2(\RR)}^2\right).
\end{split}
\end{align}
For the second term in~\eqref{e:energy}, we apply~\eqref{e:Taylor2} to bound
\begin{align*}
U\left(\absolute{1 + z_1(y)}^2\right) = U\big(\re^{2\rho(y)}\big) 
\geq \left(2 a^2 f'(a^2) - C_* \delta_2\right) |\rho(y)|^2
\end{align*}
for $y \in \RR$, yielding
\begin{align} \label{e:lowerbound2}
\int_\RR U\left(\absolute{1 + z_1(y)}^2\right)  \dd y \geq \left(2a^2 f'(a^2) - C_* \delta_2\right) \norm{\rho}_{L^2(\RR)}^2.
\end{align}
For the third term, we recall~\eqref{e:computation} and establish
\begin{align*}
&\Imag\left((1 + z_1(y)) \overline{z_1'(y)}\right) - \theta_\infty'(y) \absolute{1 + z_1(y)}^2 + \theta_\infty'(y)\\ 
&\qquad = -\theta'(y) - 2\rho(y)\left(\theta'(y) + \theta_\infty'(y)\right) -  \left(\re^{2\rho(y)} - 1 - 2\rho(y)\right)\left(\theta_\infty'(y) + \theta'(y)\right)
\end{align*}
for $y \in \RR$. Combining the latter with the bound~\eqref{e:Taylor2} and using Young's inequality with parameter $\gamma > 0$, we obtain
\begin{align} \label{e:lowerbound3}
\begin{split}
&\Absolute{\int_\RR \Imag\left((1 + z_1(y)) \overline{z_1'(y)}\right) - \theta_\infty'(y) \absolute{1 + z_1(y)}^2 + \theta_\infty'(y) \dd y} \\
&\qquad \leq \left(1 + \frac12 C_*\delta_2\right)\left(\gamma \norm{\rho}_{L^2(\RR)}^2 + \frac{1}{\gamma} \norm{\theta'+ \theta_\infty'}_{L^2(\RR)}^2\right).
\end{split}
\end{align}
Gathering the bounds~\eqref{e:first_lower},~\eqref{e:lowerbound2}, and~\eqref{e:lowerbound3}, we finally establish
\begin{align}
\label{e:bound-below-3}
\begin{split}
2E(z_1,z_2) &\geq (1-c^2) \norm{z_2}_{L^2(\RR)}^2 
+ \left(1 - c^2\right) (1 - C_* \delta_2) \left(\norm{\rho'}_{L^2(\RR)}^2 + \norm{\theta'+\theta_\infty'}_{L^2(\RR)}^2\right)\\
&\qquad + \, \left(2a^2 f'(a^2) - C_* \delta_2\right) \norm{\rho}_{L^2(\RR)}^2 \\
&\qquad - \, \absolute{k + c\omega} \left(2 + C_*\delta_2\right) \bigg(\gamma \norm{\rho}_{L^2(\RR)}^2 + \frac{1}{\gamma} \norm{\theta'+\theta_\infty'}_{L^2(\RR)}^2\bigg),
\end{split}
\end{align}
where $C_0(\delta_2) > 0$ is a constant depending on $C_*, c$ and $\delta_2$ only. By Young's inequality, we can estimate
\begin{align} \label{e:final_young}
    \norm{\theta'+\theta_\infty'}_{L^2(\RR)}^2 \geq \frac12 \|\theta'\|_{L^2(\RR)}^2 - 3\|\theta_\infty'\|_{L^2(\RR)}^2.
\end{align}

Our next step is to choose $c \in (-1,1)$, $\gamma > 0$, and $\delta_2 \in (0,1)$ in such a way that all constants in front of the norms on right-hand side of~\eqref{e:bound-below-3} are strictly positive. We distinguish between the cases $f(a^2) > 0$, $f(a^2) < 0$, and $f(a^2) = 0$.

First, we consider the case $f(a^2) > 0$. Here,~\eqref{e:disp-rel} allows us to choose $c= -\frac{k}\omega \in (-1, 1)$ such that $k + c\omega = 0$. Subsequently choosing $\gamma = 1$ and taking $\delta_2 > 0$ sufficiently small, we find that all constants in front of the norms on the right-hand side of~\eqref{e:bound-below-3} are strictly positive by~\eqref{e:condition}. Combining the latter with~\eqref{e:final_young}, settles the lower bound~\eqref{e:energy-bound-below} in the case $f(a^2) > 0$. 

Next, suppose $f(a^2) < 0$. Here,~\eqref{e:disp-rel} and~\eqref{e:condition} allow us to choose 
\begin{align*}
c = -\frac{\omega}{k} \in (-1, 1), \qquad \gamma = -|k| \frac{a^2 f'\big(a^2\big) - 2 f\big(a^2\big)}{2 f\big(a^2\big)} > 0
\end{align*}
such that
\begin{align*}
|k+c\omega| = -\frac{f(a^2)}{|k|}, \qquad 1 - c^2 - \frac{2|k+c\omega|}{\gamma} &= \left(1-c^2\right) \frac{a^2 f'\big(a^2\big) + 2f\big(a^2\big)}{a^2 f'\big(a^2\big) - 2 f\big(a^2\big)} > 0,
\end{align*}
and
\begin{align*}
2a^2f'\big(a^2\big) - 2|k+c\omega|\gamma &= a^2f'\big(a^2\big) + 2f\big(a^2\big) > 0.
\end{align*}
Subsequently taking $\delta_2 > 0$ sufficiently small, we again obtain that all constants in front of the norms on the right-hand side of~\eqref{e:bound-below-3} are strictly positive, which, together with~\eqref{e:final_young}, yields the bound~\eqref{e:energy-bound-below} for the case $f(a^2) < 0$. 

Finally, suppose $f(a^2) = 0$. We have $\frac{k}{\omega} \in \{\pm 1\}$ by~\eqref{e:disp-rel}. Setting $c = -\frac{k}{\omega}(1 - \delta_2) \in (-1,1)$ and $\gamma = 1/\sqrt{\delta_2}$, we find
\begin{align*}
|k+c\omega| = |k| \delta_2
\end{align*}
and
\begin{align*}
2a^2 f'(a^2) - 2|k+c\omega| \gamma &= 2a^2 f'(a^2) - 2\sqrt{\delta_2} \, |k|, \qquad 1-c^2-\frac{2|k+c\omega|}{\gamma} = 2\delta_2 - 2\delta_2 \sqrt{\delta_2} \, |k| - \delta_2^2.
\end{align*}
Hence, taking $\delta_2 > 0$ sufficiently small and recalling~\eqref{e:condition}, we find that all constants in front of the norms on the right-hand side of~\eqref{e:bound-below-3} are strictly positive. Combining this with~\eqref{e:final_young}, establishes~\eqref{e:energy-bound-below} for the case $f(a^2) = 0$, which concludes the proof.
\end{proof}

\section{Proof of the main result} \label{sec:mainresultproof}

We prove our orbital stability result, Theorem~\ref{t:main}. The proof relies on the conservation of the energy functional~\eqref{e:energy} and employs the upper and lower bounds obtained in Lemmas~\ref{l:finite-energy} and~\ref{l:energy-controls-norm}. 

\begin{proof}[Proof -- Theorem~\ref{t:main}.] Let $c \in (-1,1)$ be as in Lemma~\ref{l:energy-controls-norm}, and let $\delta > 0$. Suppose $\theta_\infty \in \cC^1(\RR,\RR)$, $w_0 \in H^2(\RR)$, and $v_0 \in H^1(\RR)$ satisfy~\eqref{e:initial-smallness}. 

We start by constructing a classical solution $u(t)$ to~\eqref{e:KG} with initial condition~\eqref{e:IC}. The existence of such a solution follows from the local well-posedness result for the perturbation equation~\eqref{e:short-time-dynamic} in $Z$-coordinates, established in~\S\ref{sec:local_existence}. Thus, recalling the coordinate transformations~\eqref{e:coord},~\eqref{e:def_zvariable}, and~\eqref{e:defZ}, we introduce associated initial data
\begin{align*}
z_1(y) = \frac{w_0(y)}{a \re^{\ri k y + \ri \theta_\infty(y)}}, \qquad z_2(y) = \frac{1}{\sqrt{1-c^2}} \left(\frac{v_0(y) - \ri \left(c k + \omega\right) w_0(y) + c w_0'(y)}{a  \re^{\ri k y + \ri \theta_\infty(y)}} + c \ri\theta_\infty'(y)\right).
\end{align*}
We note that $z_1 \in H^2(\RR)$ and $z_2 \in H^1(\RR)$. Moreover, there exists an $E_0$-independent constant $C_0 > 0$ such that
\begin{align} \label{e:boundsz1z2}
\|z_1\|_{H^1(\RR)} \leq C_0E_0 \leq C_0\delta, \qquad \|z_2\|_{L^2(\RR)} \leq C_0E_0 \leq 2 C_0 \delta.
\end{align}
Thanks to Proposition~\ref{p:short-time-existence}, there exist a maximal time $T \in (0,\infty]$ and a classical solution $Z \in \cC([0,T),X) \cap \cC^1([0,T),D(\Lambda))$ to~\eqref{e:short-time-dynamic} with initial value $Z(0) = \transp{(z_1,z_2)}$ such that, if $T < \infty$, then~\eqref{e:blow-up} holds. Let $z(t)$ be the first coordinate of $Z(t)$. Since the second coordinate of $Z(t)$ is given by the temporal derivative $z_t(t)/\sqrt{1-c^2}$, we find that $z$ obeys~\eqref{e:regz}. Thus, by construction, the function $u \colon \RR \times [0,T) \to \RR$ given by 
\begin{align} \label{e:defu}
u(x,t) = \left(1 + z(x-ct,t)\right) a\re^{\ri k x + \ri \omega t + \ri \theta_\infty(x-ct)}
\end{align}
satisfies~\eqref{e:regu} and is a classical solution to the complex Klein-Gordon equation~\eqref{e:KG} with initial condition~\eqref{e:IC}.

Our next step is to decompose $z(t)$ into polar coordinates and obtain bounds on the radius and phase in terms of $z(t)$. To do so, we apply Proposition~\ref{p:polar-coordinate}, which first requires expressing the initial condition $z(0) = z_1$ in polar coordinates. To this end, let $\delta_2 \in (0,1)$ be as in Lemma~\ref{l:energy-controls-norm} and fix $\delta_1 \in (0,\frac14 \delta_2)$. Taking $\delta > 0$ smaller if necessary and using~\eqref{e:boundsz1z2} and the continuous embedding $H^1(\RR) \hookrightarrow L^\infty(\RR)$, we arrange for $\|z_1\|_{L^\infty(\RR)} \leq \delta_1$. Since the principal logarithm $h \colon \mathbb{D}_1 \to \CC$, $h(s) = \log(1+s)$ is an analytic function on the open disk $\mathbb{D}_1$ of radius $1$ and satisfies $h(0) = 0$, the function $\log(1+z_1)$ is a well-defined element of the algebra $H^2(\RR)$ via the holomorphic functional calculus. We define $\rho_0,\theta_0 \in H^2(\RR,\RR)$ as its real and imaginary part through
\begin{align} \label{e:polardecomplog}
\rho_0 + \ri \theta_0 = \log(1 + z_1).
\end{align}
To bound the $L^2$-norms of $\rho_0$ and $\theta_0$, we use the mean value theorem, yielding a constant $C_* > 0$ such that $|\log(1+s)| \leq C_* |s|$ for $s \in \mathbb{D}_{\frac12}$. Consequently, the estimate~\eqref{e:boundsz1z2} implies
\begin{align} \label{e:ICbounds}
\norm{\rho_0}_{L^2(\RR)}, \norm{\theta_0}_{L^2(\RR)} \leq C_* \norm{z_1}_{L^2(\RR)} \leq C_* C_0 E_0.
\end{align}
On the other hand, we use the reverse triangle inequality to infer
\begin{align*}
\norm{\re^{\rho_0} - 1}_{L^\infty(\RR)} = \norm{\Absolute{z_1+1} - 1}_{L^\infty(\RR)} \leq \norm{z_1}_{L^\infty(\RR)} \leq \delta_1.
\end{align*}
Hence, an application of Proposition~\ref{p:polar-coordinate} yields an $E_0$-independent constant $C_1 > 0$ and unique continuously differentiable functions $\rho,\theta \colon \RR \times [0,t_{\max}) \to \RR$ with $\rho(0) = \rho_0$ and $\theta(0) = \theta_0$ such that the polar decomposition
\begin{align} \label{e:polar_proof}
z(y,t) = \re^{\rho(y,t) + \ri \theta(y,t)} - 1 
\end{align}
and the bounds
\begin{align}
\label{e:prop_bounds}
\begin{split}
\norm{\rho(t)}_{L^\infty(\RR)} &\leq 4\delta_1 \leq \delta_2,\\
\norm{\rho_y(t)}_{L^2(\RR)}^2 + \norm{\theta_y(t)}_{L^2(\RR)}^2 &\leq C_1^2 \norm{z_y(t)}_{L^2(\RR)}^2,\\
\norm{\rho_t(t)}_{L^2(\RR)}^2 + \norm{\theta_t(t)}_{L^2(\RR)}^2 &\leq C_1^2 \norm{z_t(t)}_{L^2(\RR)}^2,\\
\norm{\theta(t)}_{L^2(\RR)} &\leq \norm{\theta_0}_{L^2(\RR)} + C_1 t  \sup_{s \in [0,t]} \norm{\partial_s z(s)}_{L^2(\RR)}
\end{split}
\end{align}
hold for all $t\in [0, t_{\max})$, where we denote
\begin{align*}
t_{\max} = \sup \Set{t \in [0, T) : \big\lVert\, \absolute{1 + z(t)} - 1\big\rVert_{L^\infty(\RR)} \leq 2 \delta_1}.
\end{align*}
In particular,~\eqref{e:boundsz1z2} and~\eqref{e:prop_bounds} imply
\begin{align} \label{e:ICbounds2}
\norm{\rho_0'}_{L^2(\RR)}, \norm{\theta_0'}_{L^2(\RR)} \leq C_1 C_0 E_0.
\end{align}

We proceed with bounding the norms of $\rho(t)$, $\theta_y(t)$, and $z_t(t)$ for $t \in [0,t_{\mathrm{\max}})$ in terms of their initial data, using the energy as a Lyapunov functional. We note that, by our choice of $c$, we have $k+c\omega = 0$ whenever $f(a^2) > 0$. Moreover, from~\eqref{e:prop_bounds}, we know that $\norm{\rho(t)}_{L^\infty(\RR)} \leq \delta_2 \leq 1$ for all $t \in [0,t_{\mathrm{\max}})$. This ensures that the energy $E(Z(t))$ is well-defined for all $t \in [0,t_{\max})$ by Lemma~\ref{l:finite-energy}. So, Lemma~\ref{l:conservation-energy} yields
\begin{align} \label{e:energy_eq}
E(Z(t)) = E(Z_0)
\end{align} 
for all $t\in [0, t_{\max})$. The polar decompositions~\eqref{e:polardecomplog} and~\eqref{e:polar_proof} along with the first bound in~\eqref{e:prop_bounds} allow us to apply Lemmas~\ref{l:finite-energy} and~\ref{l:energy-controls-norm} to bound the right-hand side of~\eqref{e:energy_eq} from above and the left-hand side of~\eqref{e:energy_eq} from below, respectively. Recalling that the second coordinate of $Z(t)$ is given by $z_t(t)/\sqrt{1-c^2}$ and using~\eqref{e:boundsz1z2},~\eqref{e:ICbounds}, and~\eqref{e:ICbounds2}, this results in $E_0$-independent constants $C_2,C_3 > 0$ such that
\begin{align} \label{e:key}
\begin{split}
&\norm{z_t(t)}_{L^2(\RR)}^2 + \norm{\rho(t)}_{H^1(\RR)}^2 + \norm{\theta_y(t)}_{L^2(\RR)}^2 \\
&\qquad \leq  C_2 \bigg(\norm{z_2}_{L^2(\RR)}^2 + \norm{\rho_0}_{H^1(\RR)}^2 + \norm{\theta_0'}_{L^2(\RR)}^2 + \norm{\theta_\infty'}_{L^2(\RR)}^2\bigg) \leq C_3 E_0^2 \leq 2C_3 \delta^2
\end{split}
\end{align}
for $t \in [0,t_{\max})$. 

Next, we argue that the solution $Z(t)$ is global, i.e., that it holds $t_{\max} = T = \infty$. First, we note that the mean value theorem yields a constant $C_4 > 0$ such that
\begin{align} \label{e:taylor4}
\Absolute{\re^s - 1} \leq C_4|s|, \qquad \Absolute{\re^s} \leq C_4
\end{align}
for $s \in \mathbb{D}_1$. Hence, using the continuous embedding $H^1(\RR) \hookrightarrow L^\infty(\RR)$ and taking $\delta > 0$ smaller if necessary, estimate~\eqref{e:key} implies
\begin{align} \label{e:zest}
\begin{split}
\norm{\, \Absolute{1+z(t)}-1}_{L^p(\RR)} &= \norm{\re^{\rho(t)} - 1}_{L^p(\RR)} \leq \delta_1, \\ \norm{z_y(t)}_{L^2(\RR)} &\leq \norm{\re^{\rho(t)}}_{L^\infty(\RR)}\left(\norm{\rho_y(t)}_{L^2(\RR)} + \norm{\theta_y(t)}_{L^2(\RR)}\right) \leq \delta_1
\end{split}
\end{align}
for all $t \in [0,t_{\max})$ and $p \in \{2,\infty\}$. Combining the latter with the fact that $z \colon [0,T) \to H^1(\RR) \hookrightarrow L^\infty(\RR)$ is continuous by~\eqref{e:regz}, it follows that, by definition, $t_{\max}$ must be equal to $T$. Consequently, the inequalities~\eqref{e:key} and~\eqref{e:zest} preclude~\eqref{e:blow-up} and it must hold $\infty = T = t_{\max}$.

We proceed with establishing the decomposition~\eqref{polardecomp} and the bounds~\eqref{e:mainbounds}. Combining~\eqref{e:defu} and~\eqref{e:polar_proof}, we find that $u(x,t)$ satisfies~\eqref{polardecomp}, where $\varrho,\vartheta \colon \RR \times [0,\infty) \to \RR$ are defined by
\begin{align*}
\varrho(x,t) = \rho(x-ct,t)
\end{align*}
and
\begin{align*}
\vartheta(x,t) = \theta(x-ct,t) - \theta_\infty(x) + \theta_\infty(x-ct)
= \theta(x-ct,t) - \int_{-ct}^0 \theta_\infty'(x + y) \dd y.
\end{align*}
Clearly, $\varrho$ and $\vartheta$ are continuously differentiable in space and time, since $\rho$, $\theta$, and $\theta_\infty$ are also $\cC^1$. On the one hand, the identities
\begin{align} \label{e:rhoids}
\begin{split}
\norm{\varrho(t)}_{H^1(\RR)} = \norm{\rho(t)}_{H^1(\RR)}, \qquad
\norm{\varrho_t(t)}_{L^2(\RR)} \leq \norm{\rho_t(t)}_{L^2(\RR)} + |c| \, \norm{\rho_y(t)}_{L^2(\RR)} 
\end{split}
\end{align}
show that $\varrho \in H^1(\RR)$ and $\varrho_t \in L^2(\RR)$ for $t \geq 0$. On the other hand, the fact that $\vartheta \in H^1(\RR)$ and $\vartheta_t \in L^2(\RR)$ follows from the estimates 
\begin{align} \label{e:thetaids}
\begin{split}
\norm{\vartheta(t)}_{L^2(\RR)} &\leq \norm{\theta(t)}_{L^2(\RR)} +\left(\int_{-|c|t}^0 \int_{-|c|t}^0 \int_{\RR} |\theta_\infty'(x+y) \theta_\infty'(x+l)| \, \dd x \dd y \dd l\right)^{\frac12} \\
&\leq \norm{\theta(t)}_{L^2(\RR)} + |c| t \, \norm{\theta_\infty'}_{L^2(\RR)}
\end{split}
\end{align}
and
\begin{align} \label{e:thetaids2}
\begin{split}
\norm{\vartheta_x(t)}_{L^2(\RR)} &\leq \norm{\theta_y(t)}_{L^2(\RR)} + 2 \norm{\theta_\infty'}_{L^2(\RR)},\\ 
\norm{\vartheta_t(t)}_{L^2(\RR)} &\leq \norm{\theta_t(t)}_{L^2(\RR)} + |c| \left(\norm{\theta_y(t)}_{L^2(\RR)} + \norm{\theta_\infty'}_{L^2(\RR)}\right)
\end{split}
\end{align}
for $t \geq 0$. Finally, combining~\eqref{e:rhoids},~\eqref{e:thetaids}, and~\eqref{e:thetaids2} with~\eqref{e:ICbounds},~\eqref{e:key}, and the last two inequalities in~\eqref{e:prop_bounds} yields~\eqref{e:mainbounds}. 

It remains to show~\eqref{e:orbitalstab}. To this end, take $R>0$, $t \geq 0$, and $x_* \in \RR$. We set $\psi(x,t) := \vartheta(x, t) + \theta_\infty(x)$, $\gamma = \psi(x_*, t)$, and $I = [x_*-R,x_*+R]$. Using~\eqref{e:key},~\eqref{e:thetaids2}, and
\begin{align*}
\psi(x,t) - \gamma = \int_{x_*}^x \psi_x(y,t) \, \dd y, \qquad x \in \RR,
\end{align*}
we find an $E_0$-, $t$-, $x_*$- and $R$-independent constant $C_5 > 0$ such that
\begin{align} \label{e:psiest}
\begin{split}
\norm{\psi(t)  - \gamma}_{L^2(I)} &\leq \left(\int_I \int_I \int_I |\psi_x(y,t)| |\psi_x(l,t)| \, \dd y \dd l \dd x\right)^{\frac12}\\ &\leq 2\sqrt{R} \, \norm{\psi_x(t)}_{L^2(I)} \leq C_5 \sqrt{R} \, E_0 \leq 2 C_5 \sqrt{R} \, \delta.
\end{split}
\end{align}
Hence, taking $\delta > 0$ smaller if necessary, applying~\eqref{e:key},~\eqref{e:taylor4},~\eqref{e:rhoids},~\eqref{e:thetaids2}, and~\eqref{e:psiest}, and using the embedding $H^1(I) \hookrightarrow L^\infty(I)$, we infer
\begin{align*}
\Absolute{u(x,t) - a \re^{\ri k x + \ri \omega t + \ri \gamma}} &= a \Absolute{\re^{\varrho(x, t) + \ri \psi(x, t)  - \ri \gamma} - 1} \leq a C_4 \Absolute{\varrho(x, t) + \ri \psi(x, t)  - \ri \gamma},\\
\Absolute{u_x(x,t) - a \ri k \re^{\ri k x + \ri \omega t + \ri \gamma}} &\leq |k| \Absolute{u(x,t) - a \re^{\ri k x + \ri \omega t + \ri \gamma}} + C_4\left(\Absolute{\varrho_x(x,t)} + \Absolute{\psi_x(x,t)}\right)
\end{align*}
for $x \in I$. Taking $L^2$-norms of the latter two inequalities and applying~\eqref{e:mainbounds} and~\eqref{e:psiest} yields~\eqref{e:orbitalstab}, which completes the proof.
\end{proof}

\section{Numerics} \label{sec:numerics}

\begin{figure}[b!]
\centering
\includegraphics[width = \textwidth]{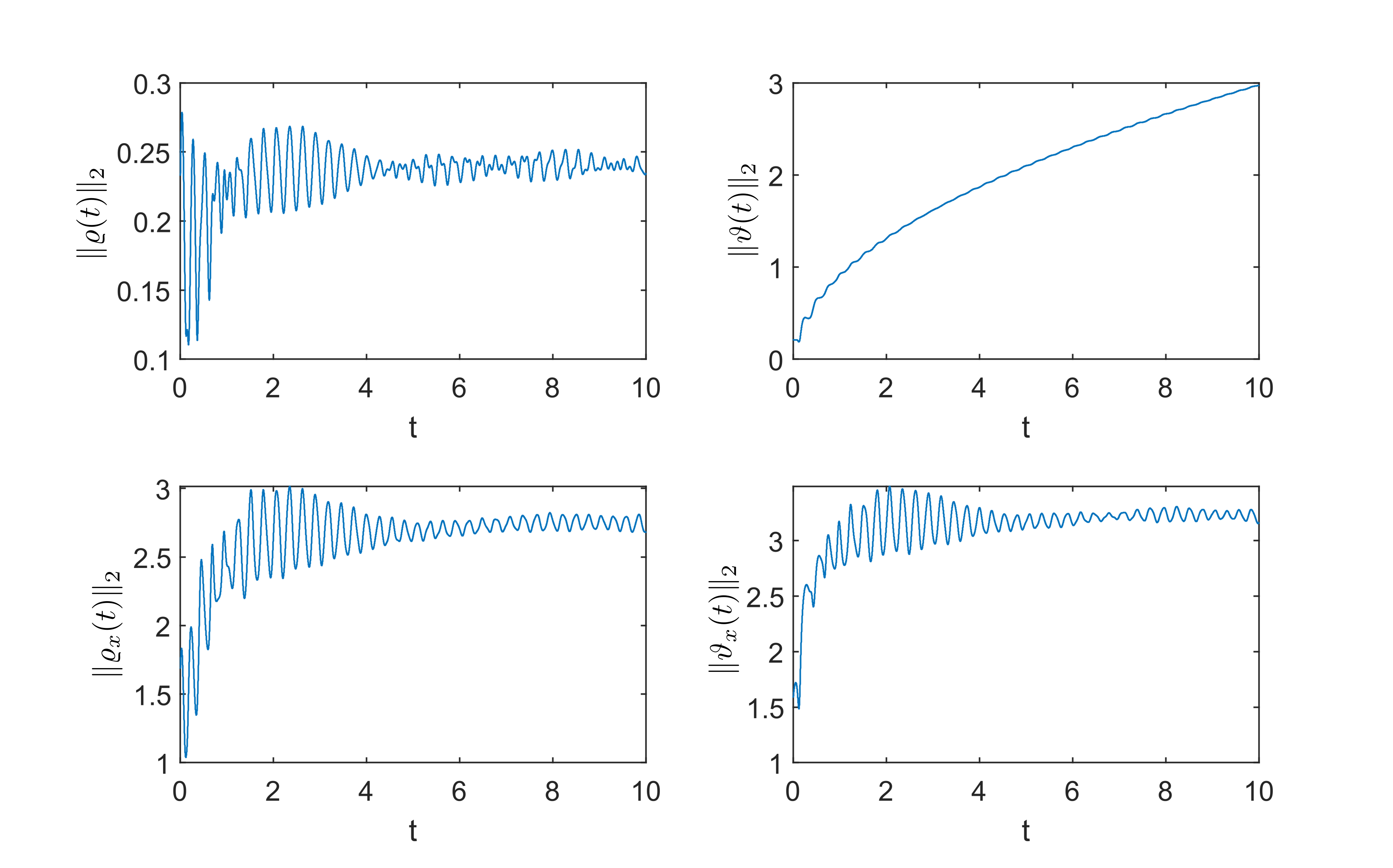}
\caption{Plots of the $L^2$-norms of the polar coordinates $\varrho(t)$ and $\vartheta(t)$ of the perturbed plane-wave solution~\eqref{e:planewavepert} with parameters~\eqref{e5.1} and initial condition~\eqref{e:IC}, along with those of their spatial derivatives. The norms of $\varrho(t)$, $\varrho_x(t)$ and $\vartheta_x(t)$ remain bounded in time. In contrast, the top-right panel illustrates that the $L^2$-norm of $\vartheta(t)$ grows over time.} \label{fig1}
\end{figure}

We discuss the outcome of numerical simulations pertaining to our main result, Theorem~\ref{t:main}. Here, we focus on the \emph{defocusing cubic} complex Klein-Gordon equation, i.e., we restrict ourselves to the case $f(s)=1+s$. All computations are carried out in MATLAB. We consider plane-wave solutions~\eqref{e:plane-wave} which are perturbed by complex-valued Gaussians. Specifically, we consider a perturbed plane-wave solution 
\begin{align} \label{e:planewavepert}
u(x,t) = a \re^{\ri k x + \ri \omega t} \, \re^{\varrho(x,t) + \ri \vartheta(x,t)}, 
\end{align}
to~\eqref{e:KG} with initial condition~\eqref{e:IC} with parameters
\begin{align} \label{e5.1}
\begin{split}
\omega=10,&\qquad k=2\pi, \qquad a^2 = \omega^2 - k^2 - 1,\\
\theta_\infty(x) = 0, \qquad  w_0(x) &= 4(1+\ri)\eu^{-25(x-10)^2}, \qquad v_0(x) = 40(1+\ri)\eu^{-25(x-10)^2}.
\end{split}
\end{align}
We emphasize that the existence condition~\eqref{e:disp-rel} and the spectral condition~\eqref{e:condition} are satisfied. To compute the time evolution of the perturbed plane-wave solution~\eqref{e:planewavepert}, we employ a Strang splitting algorithm. We simulate the dynamics on a spatial domain of $20$ units with periodic boundary conditions. The simulation time is chosen such that the perturbation does not reach the domain boundary, preventing contributions of boundary effects. The results are plotted in~\hyperref[fig1]{Fig.\ref*{fig1}}.

\begin{figure}[t!]
\centering
\includegraphics[width = 0.6\textwidth]{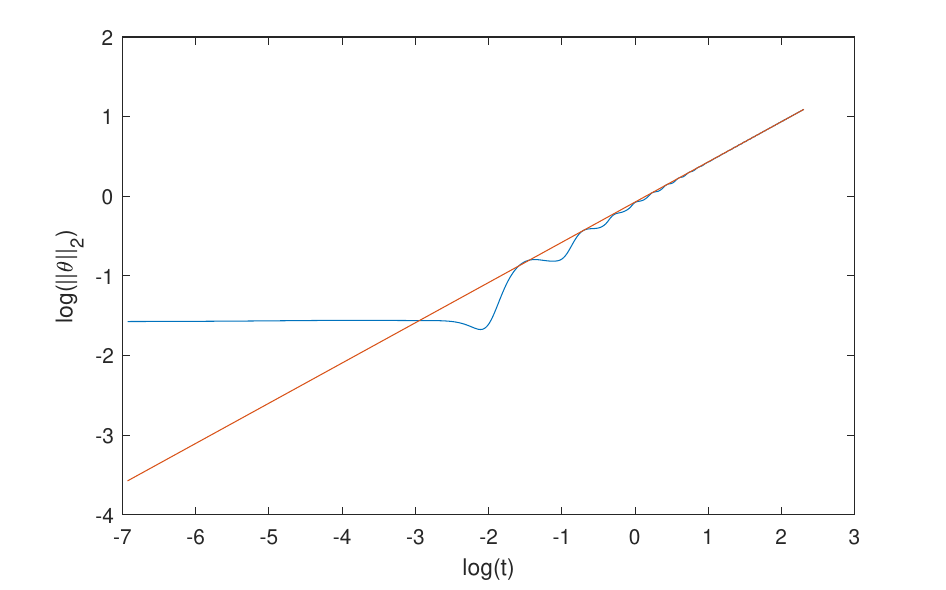}
\caption{Log-log plot of the $L^2$-norm of the phase $\vartheta(t)$ of the perturbed plane-wave solution~\eqref{e:planewavepert} with parameters~\eqref{e5.1} and initial condition~\eqref{e:IC}. The best linear fit, shown in red, is computed using the last quarter of data points. Its slope is approximately $0.5051$, suggesting that $\norm{\vartheta(t)}_{L^2(\mathbb{R})}$ grows with rate $\sqrt{t}$.} \label{fig2}
\end{figure}

\begin{figure}[t!]
\centering
\includegraphics[width = 0.95\textwidth]{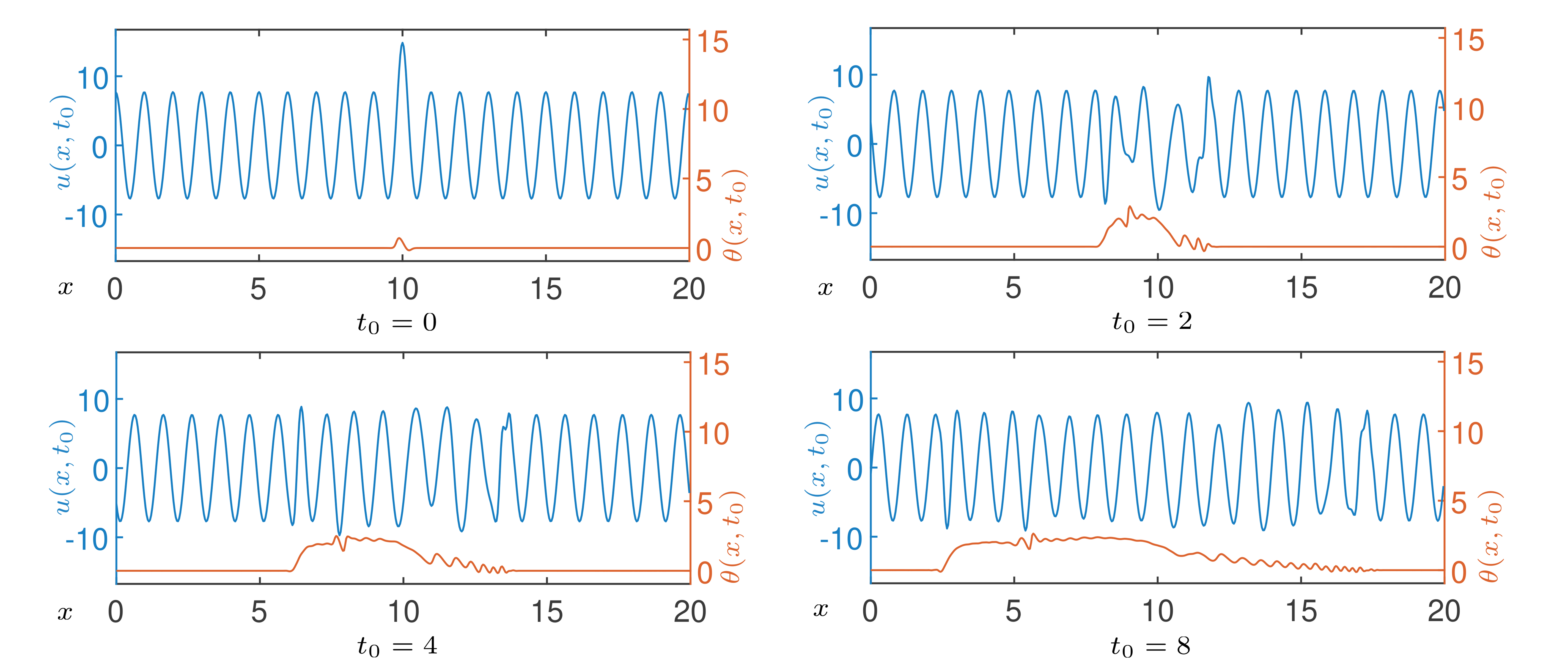}
\caption{Plots of the real part of the perturbed plane-wave solution~\eqref{e:planewavepert} (blue), with parameters~\eqref{e5.1} and initial condition~\eqref{e:IC}, along with its phase $\vartheta(t_0)$ (orange) at different times $t_0$. The initial perturbations are given by~\eqref{e5.1} scaled by $\smash{\frac{7}{4}}$ for improved visibility. The perturbation triggers a phase defect that travels outwards in both spatial directions.} \label{fig3}
\end{figure}

The numerical computations indicate that the bounds~\eqref{e:mainbounds} on the $L^2$-norm of $\varrho(t)$, $\varrho_x(t)$, and $\vartheta_x(t)$ in~\eqref{e:mainbounds} are sharp, as these do not seem to decay over time. However, the numerics also suggest that the linear growth bound for $\smash{\norm{\vartheta(t)}_{L^2(\RR)}}$ in~\eqref{e:mainbounds} is not optimal. Indeed, a linear regression for the log-log plot of $\|\vartheta(t)\|_{L^2(\RR)}$ indicates that $\|\vartheta(t)\|_{L^2(\RR)}$ grows algebraically with rate $\smash{\sqrt{t}}$, see~\hyperref[fig2]{Fig.~\ref*{fig2}}. The $L^\infty$-norm of $\vartheta(t)$, on the other hand, remains bounded. The dynamics of the perturbed plane-wave solution $u(t)$ together with its phase $\vartheta(t)$ are illustrated in~\hyperref[fig3]{Fig.~\ref*{fig3}}. Similar behavior was observed in numerical simulations for different values of $\omega$, $k$, and initial perturbations $w_0, v_0$.

\section{Spectral stability analysis} \label{sec:specstab}

The analysis in this section confirms that our orbital stability result, Theorem~\ref{t:main}, is sharp in the sense that it holds up to the spectral stability boundary. To study the spectral stability of the plane-wave solution~\eqref{e:plane-wave}, we consider the perturbation
\begin{align} \label{e:pert-plane-wave}
u(x,t) = a\re^{\ri kx + \ri \omega t}\left(1+z(x-ct,t)\right),
\end{align}
where the free parameter $c \in \RR$ denotes the speed of the co-moving frame in which we observe the dynamics. Inserting~\eqref{e:pert-plane-wave} in the complex Klein-Gordon equation~\eqref{e:KG} while using~\eqref{e:disp-rel}, we arrive at
\begin{align} \label{e:zdynamics}
z_{tt} - \left(1 - c^2\right) z_{yy} - 2c z_{yt} + 2 \ri \omega z_t - 2 \ri\left(c\omega + k\right) z_y + \left(f\big(a^2|1+z|^2\big) - f\big(a^2\big)\right) (1+z) = 0.
\end{align}
Writing~\eqref{e:zdynamics} as a real system in $Z = (\Real z,\Imag z)^\top$ and linearizing, we obtain 
\begin{align*}
Z_{tt} + \mathcal{J}_c(\partial_y) Z_t + \mathcal{H}_c(\partial_y) Z = 0
\end{align*}
with
\begin{align*}
\mathcal{J}_c(D) = \begin{pmatrix} -2cD & -2\omega \\ 2\omega & -2cD\end{pmatrix}, \qquad \mathcal{H}_c(D) = \begin{pmatrix} -\big(1-c^2\big) D^2 + 2a^2f'\big(a^2\big) & 2(c\omega + k)D \\ -2(c\omega + k)D & -\big(1-c^2\big) D^2\end{pmatrix}.
\end{align*}

The spectrum associated with the plane-wave solution~\eqref{e:plane-wave} to~\eqref{e:KG} is thus given by the set 
\begin{align*}
\sigma_c = \left\{\lambda \in \CC : \mathcal{P}_c(\lambda) \text{ is not invertible}\right\},
\end{align*}
where $\mathcal{P}_c(\lambda) \colon H^2(\RR) \subset L^2(\RR) \to L^2(\RR)$ is the quadratic operator pencil
\begin{align*}
\mathcal{P}_c(\lambda) = \lambda^2 I_2 + \lambda \mathcal{J}_c(\partial_y) + \mathcal{H}_c(\partial_y).
\end{align*}
We say that the plane wave~\eqref{e:plane-wave} is \emph{spectrally stable} if $\sigma_c$ is confined to the imaginary axis and \emph{spectrally unstable} if it is not spectrally stable. Since $\mathcal{P}_c(\lambda)$ has constant coefficients, we find that it is invertible if and only if its Fourier symbol
\begin{align*}
\hat{\mathcal{P}}_c(\lambda,\ell) = \lambda^2 I_2 + \lambda \mathcal{J}_c(\ri \ell) + \mathcal{H}_c(\ri\ell)
\end{align*}
is invertible for each $\ell \in \RR$. That is, $\mathcal{P}_c(\lambda)$ is not invertible if and only if there exists $\ell \in \RR$ such that $(\lambda,\ell)$ is a root of the determinantal function $\mathcal{E}_c \colon \CC \times \RR \to \CC$ given by
\begin{align*}
\mathcal{E}_c(\lambda,\ell) = \det \hat{\mathcal{P}}_c(\lambda,\ell).
\end{align*}
Since we have $\hat{\mathcal{P}}_c(\lambda,\ell) = \hat{\mathcal{P}}_0(\lambda - c\ri \ell,\ell)$ for all $\lambda \in \CC$ and $c,\ell \in \RR$, the spectral stability of the plane wave is independent of the choice of $c \in \RR$. 

The following result, summarizing the outcome of our spectral analysis, shows that condition~\eqref{e:condition} is sharp in the sense that it coincides with a change in spectral stability. Its proof employs the instability index theory for quadratic operator pencils as developed in~\cite{Bronski2014Instability}. 

\begin{proposition}
Let $(a,k,\omega) \in (0,\infty) \times \RR^2$ be such that the existence condition~\eqref{e:disp-rel} is satisfied.
\begin{itemize}
\item For $f(a^2) > 0$ the plane-wave solution~\eqref{e:plane-wave} to~\eqref{e:KG} is spectrally stable if $f'(a^2) > 0$ and spectrally unstable if
\begin{align} \label{e:unstabcond}
-2f(a^2) < a^2 f'\big(a^2\big) < 0. 
\end{align}
\item For $f(a^2) < 0$ the plane-wave solution~\eqref{e:plane-wave} to~\eqref{e:KG} is spectrally stable if $a^2 f'(a^2) > -2f(a^2)$ and spectrally unstable if
\begin{align}\label{e:unstabcond2}
0 < a^2 f'\big(a^2\big) < -2f(a^2).
\end{align}
\end{itemize}
\end{proposition}
\begin{proof}
Assume that condition~\eqref{e:condition} holds. We show that the plane wave is spectrally stable. As in the orbital stability analysis, we choose $c \in \RR$ depending on the sign of $f(a^2)$. First, we consider the case $f(a^2) > 0$. Using~\eqref{e:disp-rel} and setting $\smash{c = -\frac{k}{\omega}} \in (-1,1)$, we find that the diagonal self-adjoint operator $\mathcal{H}_c(\partial_y)$ is positive semidefinite on $\smash{L^2_{\mathrm{per}}(0,L)}$ for any $L \geq 0$. Therefore, the instability index count in~\cite{Bronski2014Instability} implies that the plane wave is spectrally stable. Next, we consider the case $\smash{f(a^2)} < 0$. We use~\eqref{e:disp-rel} to set $\smash{c = -\frac{\omega}{k}} \in (-1,1)$. By~\eqref{e:condition} we find that
\begin{align*}
\mathcal{H}_{-\frac{\omega}{k}}(\ri k \ell) = \begin{pmatrix} \big(k^2 -\omega^2\big)\ell^2 + 2a^2 f'\big(a^2\big) & 2(k^2 - \omega^2)\ri \ell \\ -2(k^2-\omega^2)\ri \ell & \big(k^2 -\omega^2\big)\ell^2\end{pmatrix}
\end{align*}
possesses nonnegative determinant
\begin{align*}
\ell^2 \left(k^2 - \omega^2\right)\left(\ell^2\left(k^2 - \omega^2\right) + 2 a^2 f'(a^2) - 4\left(k^2 - \omega^2\right)\right)
\end{align*}
and one nonnegative eigenvalue
\begin{align*}
a^2f'(a^2) + \ell^2 \left(k^2 - \omega^2\right) + \sqrt{\left(a^2f'(a^2)\right)^2 + 4 \ell^2 (k^2 - \omega^2)^2}
\end{align*}
for all $\ell \in \RR$. Therefore, both eigenvalues of $\mathcal{H}_c(\ri k \ell)$ are nonnegative and the self-adjoint operator $\mathcal{H}_c(\partial_y)$ is positive semidefinite on $L^2_{\mathrm{per}}(0,L)$ for any $L \geq 0$. Therefore, the plane wave is spectrally stable by~\cite{Bronski2014Instability}. 

Next, we show that~\eqref{e:unstabcond} and~\eqref{e:unstabcond2} yield spectral instability. We use that $\mathcal{E}_0(\ri \omega,\ell)$ is a quartic polynomial in $\omega$ with real coefficients for $\ell \in \RR$. Let $\Delta(\ell)$ be the discriminant of $\mathcal{E}_0(\ri \cdot ,\ell)$. We compute
\begin{align*}
\Delta(0) = 0, \qquad \Delta'(0) = 0, \qquad \Delta''(0) = 1024 a^2 f'(a^2) \left(2 \omega^2 + a^2 f'\big(a^2\big)\right)^3 \left(2 \omega^2 - 2 k^2 + a^2 f'\big(a^2\big)\right).
\end{align*}
If we have~\eqref{e:unstabcond} or~\eqref{e:unstabcond2}, then $\Delta''(0)$ is negative by~\eqref{e:disp-rel}. Hence, there exists $\ell_0 > 0$ such that the discriminant $\Delta(\ell)$ is negative for $\ell \in (-\ell_0,\ell_0) \setminus \{0\}$. Therefore, the quartic polynomial $\mathcal{E}_0(\ri \cdot ,\ell)$ must possess two non-real roots for $\ell \in (-\ell_0,\ell_0) \setminus \{0\}$, implying that the plane wave is spectrally unstable. 
\end{proof}

\begin{remark}
Our analysis shows that, in the cases~\eqref{e:unstabcond} and~\eqref{e:unstabcond2}, the plane wave is spectrally unstable due to the fact that $\mathcal{E}_c(\ri \cdot,\ell)$ possesses nonreal roots for all Fourier frequencies $\ell \in (-\ell_0,\ell_0) \setminus \{0\}$ with $\ell_0 > 0$ sufficiently small. That is, the plane wave undergoes a \emph{sideband, Benjamin-Feir or long wavelength instability} as it crosses the stability boundary $a^2 f'(a^2) = 2\max\{0,-f(a^2)\}$. 
\end{remark}

\bibliographystyle{alphaabbr}
\bibliography{references}

\end{document}